\documentclass[preprint,10pt]{elsarticle}

\usepackage{graphicx}
\graphicspath{{Figures/}}
\usepackage{algorithm2e}
\RestyleAlgo{ruled}
\usepackage{natbib}
\usepackage{amsmath,amsfonts,amssymb}
\usepackage{mathrsfs}
\usepackage{lineno}
\usepackage{bm}
\usepackage{stmaryrd}
\usepackage{nicefrac}

\usepackage{scalerel}
\DeclareMathOperator*{\A}{\scalerel*{A}{\sum}}

\usepackage{caption}
\usepackage{subcaption}
\usepackage{xcolor}
\usepackage{epstopdf}

\usepackage{hyperref} 

\journal{arXiv}
\begin{document}

\begin{frontmatter}
\title{A dynamic implicit 3D material point-to-rigid body contact approach for large deformation analysis}

\author[1]{Robert E. Bird}
\author[1]{Giuliano Pretti}
\author[1]{William M. Coombs\corref{cor1}}
\author[1]{Charles E. Augarde}

\author[2]{Yaseen U. Sharif}
\author[2]{Michael J. Brown}

\author[3,4]{Gareth Carter}
\author[4]{Catriona Macdonald}
\author[4]{Kirstin Johnson}

\cortext[cor1]{Corresponding author: w.m.coombs@durham.ac.uk}
\address[1]{Department of Engineering, Durham University, Durham, DH1 3LE, UK}
\address[2]{School of Science and Engineering, University of Dundee, Fulton Building, Dundee, DD1 4HN, Scotland, UK}
\address[3]{Arup, Edinburgh, Scotland, UK}
\address[4]{British Geological Survey, Currie, Edinburgh, EH14 4BA, Scotland, UK}

\begin{abstract}

\noindent  Accurate and robust modelling of large deformation three dimensional contact interaction is an important area of engineering, but it is also challenging from a computational mechanics perspective.  This is particularly the case when there is significant interpenetration and evolution of the contact surfaces, such as the case of a relatively rigid body interacting with a highly deformable body.  The numerical challenges come from several non-linear sources: large deformation mechanics, history dependent  material behaviour and slip/stick frictional contact.  In this paper the Material Point Method (MPM) is adopted to represent the deformable material, combined with a discretised rigid body which provides an accurate representation of the contact surface.  The three dimensional interaction between the bodies is detected though the use of domains associated with each material point.  This provides a general and consistent representation of the extent of the deformable body without introducing boundary representation in the material point method.  The dynamic governing equations allows the trajectory of the rigid body to evolve based on the interaction with the deformable body and the governing equations are solved within an efficient implicit framework.  The performance of the method is demonstrated on a number of benchmark problems with analytical solutions.  The method is also applied to the specific case of soil-structure interaction, using geotechnical centrifuge experimental data that confirms the veracity of the proposed approach.    
\end{abstract}


\end{frontmatter}

\section{Introduction}

Modelling the interaction between deformable bodies under large deformations requires robust contact algorithms that are able to track and handle the evolving contact surfaces while enforcing the required contact constraints. In many areas of engineering, one of the bodies can be very stiff compared to the other, which permits a rigid body assumption for the stiffer material.  One such application area is soil-structure interaction, which is a key area of geotechnical engineering, particularly in offshore geotechnics where almost all interactions with the seabed involve large deformation processes. For example, the successful deployment of Offshore Renewable Energy (ORE) is underpinned by geotechnical engineers being able to understand and predict the interaction of supporting structures, moorings and other enabling works (cone penetration tests for site characterisation, ploughing for transmission cable installation, spudcans for jack up vessels, drag anchor penetration for moorings, etc.) with the seabed.  These interactions share similar characteristics in that they are truly three dimensional, involve large deformation soil-structure interaction, contact/friction, and history-dependent non-linear soil behaviour.  The trajectory of the penetrating body is also often coupled to the response of the soil, rather than being fully prescribed.   The combination of these characteristics makes modelling the physical response very challenging.  There are also several options in terms of how to represent the soil and the structure in such analyses.   


Solid mechanics-based numerical analysis in the ORE industry, and engineering in general, is dominated by the Finite Element Method (FEM). However the FEM struggles to analyse problems involving very large deformation due to issues associated with mesh distortion, remeshing and remapping of state variables (such as parameters for history-dependent materials).  There are numerous continuum-based methods that are designed to overcome this limitation, such as: the Coupled Eulerian Lagrangian Method (CELM, \cite{Noh1963}), the Arbitrary Lagrangian Eulerian (ALE, \cite{Hughes1981}) finite element method, the Particle Finite Element Method (PFEM, \cite{Idelsohn2004}), and the Material Point Method (MPM, \cite{sulsky1994particle}), amongst others.  Each of the methods have advantages and disadvantages depending on the specific problem under consideration (see \citet{Augarde2021} for a recent review).  An alternative to continuum-based approaches is to adopt the the Discrete Element Method (DEM, \cite{Cundall1979}) to represent both the structure and the soil. However, there are issues in calibrating suitable particle interaction models to represent bulk material behaviour, and analyses are very expensive compared to continuum-based approaches. This paper builds on the work of \citet{bird2024implicit} and adopts the MPM as the numerical method to represent soil behaviour due to the MPM's established track record of modelling large deformation geotechnical problems (see \citet{Solowski2021} and \citet{Vaucorbeil2020} for review articles).  The large deformation continuum formulation adopted in this paper is based on that of \citet{charlton2017igimp} and the open source Material Point Learning Environment (AMPLE) code \cite{coombs2020aample}, and is outlined in Section~2, which also explains the MPM large deformation calculations steps.  

There are also choices in terms of how to enforce contact between the soil and the structure.  The original MPM includes a form of no-slip contact in that multiple bodies interact on a shared background grid, where the governing equations are assembled and solved.  However, most of these methods struggle to rigorously enforce contact as they lack explicit definition of the contact surface between the bodies \cite{acosta2021development,pretti2024continuum}.  The focus in this paper is on soil-structure interaction where the stiffness of the structure is sufficiently large, relative to the deformable soil body, that it can be assumed to be rigid.  This allows the contact surface to be based on the geometry of the rigid body \cite{bird2024implicit}.  Therefore, this introduction outlines specific MPM contributions that focus on soil-structure interaction with a clearly defined contact boundary rather than attempt to provide a comprehensive review of all MPM-based contact formulations (see \citet{bird2024implicit} for a more detailed overview).  Within this context, \citet{nakamura2021particle} used the boundary discretisation of an infinitely rigid body to define the contact surface, while \citet{Lei2022} employed a finite element discretisation that allowed the structure to deform whilst providing an explicit representation of the contact surface.  It is worth noting that in some cases it is possible to maintain a rigid body conforming background mesh, which allows the contact constraints to be imposed directly on the background grid (see for example \citet{Martinelli2021}).  However, in all but the simplest of cases this necessitates the use of unstructured background grids, usually simplex elements, which causes other issues such as cell crossing instabilities as MPM formulations adopting $C_1$ continuous basis functions (such as generalised interpolation \cite{bardenhagen2004generalized} and B-spline \cite{Tielen2017} basis functions) are only available for structured grids.  This paper extends the two-dimensional quasi-static contact approach of \citet{bird2024implicit} to three-dimensional dynamic analysis, including coupled soil reaction dependent rigid body motion.  The deformable body is represented by Generalised Interpolation Material Points (GIMPs) in order to mitigate cell crossing instabilities.  Numerical methods for enforcing contact constraints can be broadly grouped into penalty, augmented Lagrangian and Lagrange multiplier methods \cite{wriggers2006computational}. Again following \citet{bird2024implicit}, here the geometrically non-linear contact constraints are enforced via a penalty approach and the dynamic coupled soil deformation-rigid body motion equations resolved using an implicit monolithic solver.  The contact detection and enforcement algorithm is based on the interaction between the rigid body and the generalised interpolation domain associated with each material point.  This means that the contact penalty forces are consistent with the internal forces generated based on the stress state of material points and that no additional material point-based boundary representation is required. Details of the approach are covered in Sections 3, 4 and 5 which provide a description of the rigid body, details of the contact formulation and information on specific aspects of the numerical implementation, respectively.          

The key contribution of this paper is a new three-dimensional MPM-based soil-structure interaction approach for dynamic problems that allows for the motion of the structure to evolve based on the response of the deformable body.  This is important for several areas of offshore geotechnical engineering where the motion of the rigid body is only partially prescribed.  Solving the coupled equations implicitly permits the use of large time steps compared to most of the MPM literature that focuses on explicit solutions methods. Enforcing the contact constraints based on the domains associated with each material point provides a robust, general and consistent soil-structure interaction approach. The advantages of the proposed formulation are demonstrated by a number of benchmark and exemplar test cases.    
 

\section{Material point method}

\begin{figure}[!h]
\centering
    \begin{subfigure}[t]{0.33\textwidth}
        \centering
        \includegraphics[width=0.99\textwidth]{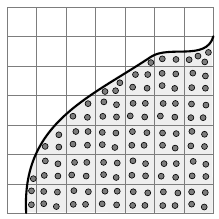}
        \caption{\footnotesize initial position}
    \end{subfigure}%
    \begin{subfigure}[t]{0.33\textwidth}
        \centering
        \includegraphics[width=0.99\textwidth]{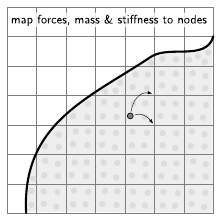}
        \caption{\footnotesize point-to-grid map}
    \end{subfigure}%
    \begin{subfigure}[t]{0.33\textwidth}
        \centering
        \includegraphics[width=0.99\textwidth]{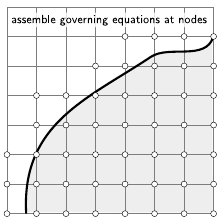}
        \caption{\footnotesize assembly}
    \end{subfigure} \newline
    \begin{subfigure}[t]{0.33\textwidth}
        \centering
        \includegraphics[width=0.99\textwidth]{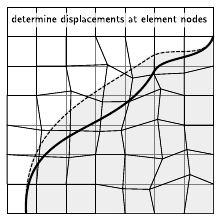}
        \caption{\footnotesize grid solution}
    \end{subfigure}%
    \begin{subfigure}[t]{0.33\textwidth}
        \centering
        \includegraphics[width=0.99\textwidth]{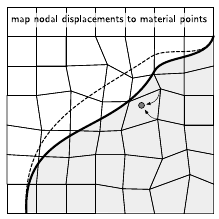}
        \caption{\footnotesize grid-to-point map}
    \end{subfigure}%
    \begin{subfigure}[t]{0.33\textwidth}
        \centering
        \includegraphics[width=0.99\textwidth]{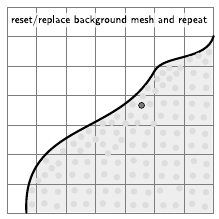}
        \caption{\footnotesize update points \& reset grid}\label{fig: method steps f}
    \end{subfigure}
\caption{Material point method steps (adapted from \citet{coombs2020aample}).}
\label{Fig:MPMsteps}
\end{figure}

In order to avoid mesh distortion issues seen in other mesh-based methods, the MPM requires two spatial discretisations: material points representing the physical body and a background grid to solve the governing equations.   A typical time step in the MPM can be broken down as shown in Figure~\ref{Fig:MPMsteps}:
\begin{enumerate}
    \item[(a)] Initial position: At the start of an analysis the physical body being analysed (shown by the light grey shaded region with a boundary defined by the black line) is discretised as a collection of material points (shown by the dark grey shaded circles) with associated volume, mass and constitutive parameters (for example Young's modulus).  The body lies within a background grid of sufficient size to cover the full extent of the material points.  
    \item[(b)] Point-to-grid mapping: Information held at material points, such as mass, stiffness and forces associated with gravitational loads, any tractions applied to material points and the stress in the material is mapped to the nodes of the background grid.  
    \item[(c)] Governing equation assembly: The governing equations describing the physical problem being analysed are assembled at the degrees of freedom of the background grid.  This leads to a finite element-like system of equations, which will change size depending on the interaction between the material points and the grid as only some nodes will be \emph{active} in the analysis (those with contributions from material points, as shown by the white shaded circles).
    \item[(d)] Grid solution: The governing equations, combined with any boundary conditions imposed directly on the grid or from immersed constraints, are solved at the \emph{active} degrees of freedom to determine the primary unknowns of the system of equations (for static and dynamic stress analysis this would usually be displacements and accelerations, respectively). How the equations are solved defines different MPMs and depends on the nature of the equations and the adopted time stepping algorithm.
    \item[(e)] Grid-to-point mapping: The solution is mapped from the grid to the material points, such as displacement, velocity, stress, deformation, volume, etc. 
    \item[(f)] Point/grid update: The positions of the material points are updated and the background grid reset or replaced.  
\end{enumerate}
These steps describe a general material point time step for quasi-static or dynamic stress analysis problems.  Details of the calculation steps will change depending on the solution procedure.  For example, for implicit implementations Step (d) will be an iterative process to find equilibrium within a given tolerance requiring several grid-to-point and point-to-grid mappings to update the deformation and stress state at material points and project the associated stiffness and internal force to the grid degrees of freedom.  For explicit implementations Step (d) will be a single step forward projection, albeit usually with a significantly smaller time step.  Explicit schemes also need to consider when to update the stress at the material points within the overall algorithm, distinguishing several methods such as Update Stress First (USF), Update Stress Last (USL) and Modified Update Stress Last (MUSL, \cite{Sulsky1995}). This paper is focused on dynamic stress analysis of solid materials, with an implicit Newmark time discretisation.  Key aspects of the formulation are provided in the following sections.  These sections are deliberately brief, with specific references for full details, as the focus on the paper is on a new three-dimensional contact approach.       

\subsection{Continuum formulation}

The continuum formulation adopted in this paper is the same as the quasi-static open-source AMPLE (A Material Point Learning Environment) code \cite{coombs2020aample}, but with the addition of inertia effects.  The updated Lagrangian formulation can be defined by the following weak statement of equilibrium 
\begin{equation}\label{eqn:weak}
  \int_{\varphi_t(\Omega)}\Bigl( \sigma_{ij}(\nabla_x \eta)_{ij} -(b_i-\rho\dot{v}_i)\eta_i \Bigr) \text{d}V - \int_{\varphi_t(\partial\Omega)}\bigl(t_i\eta_i \bigr) \text{d}s = 0.
\end{equation}
where $\dot{v}_i$ and $\varphi_t$ are the acceleration and motion of the physical body with domain, $\Omega$, which is subjected to tractions, $t_i$, on its boundary, $\partial\Omega$ (with surface, $s$), and body forces, $b_i$, acting over its volume, $V$, with density, $\rho$.  This leads to a Cauchy stress field, $\sigma_{ij}$, through the body. The weak form is derived in the current (updated, deformed) frame assuming a field of admissible virtual displacements, $\eta_i$.   
Full details of the large deformation elasto-plastic continuum mechanics formulation used in this paper can be found in \citet{charlton2017igimp} and \citet{coombs2020aample}. Within this approach, the deformation gradient is multiplicatively decomposed into elastic and plastic components and combined with  a linear relationship between Kirchhoff stress and logarithmic elastic strain along with an exponential map of the plastic flow rule.  The combination of these ingredients allows isotropic small strain plasticity algorithms to be used directly within a large deformation setting without modifying the stress update procedure \cite{Simo1992a}.  The approach is widely used in large deformation finite element methods (see for example \citet{SouzaNeto}), and in material point methods \cite{charlton2017igimp,Coombs2020continuum,coombs2020aample,cortis2018imposition,wang2021efficient,Wang2019,bird2024implicit}.    

\subsection{Material point and background grid spatial discretisation}

Introducing a background grid of finite elements, $K$, allows the weak statement of equilibrium (\ref{eqn:weak}) to be discretised as
\begin{equation}\label{eqn:MPMweak}
  \int_{\varphi_t(K)}[\nabla_x S_{v}]^{T}\{\sigma\} \text{d}V 
  -\int_{\varphi_t(K)}[S_{v}]^{T}\{b\} \text{d}V 
  +\int_{\varphi_t(K)}\rho[S_{v}]^{T}\{\dot{v}\} \text{d}V
  -\int_{\varphi_t(\partial \Omega)}\left(\{F_{N,v} ^{\partial\Omega}\} + \{F_{T,v} ^{\partial\Omega}\}\right)\text{d}s
  =\{0\},
\end{equation}
where  $[S_{v}]$ contains the basis functions that map information from the nodes (or vertices, $v$) of the background grid to other locations within the grid.  $[\nabla_x S_{v}]$ is the strain-displacement matrix containing derivatives of the basis functions with respect to the updated coordinates.  The last term is the traction comprising the normal and tangential contact forces for the deformable body, respectively $\{F_{N,v} ^{\partial\Omega}\}$ and $\{F_{T,v} ^{\partial\Omega}\}$. Note that the traction term remains integrated over the boundary of the physical body, $\partial \Omega$, in general the domain boundary will not coincide with the background mesh and traction boundary conditions cannot be applied on the edges of elements on the mesh boundary. Therefore special treatment of the application of contact boundary conditions are required to map boundary conditions from the domain's surface to the mesh boundary, this is described in Section \ref{sec: CPP}.

In the MPM the physical body being analysed is represented by a number of discrete material points with an associated volume, $V_p$, and mass, $m_p$.  This allows the volumetric integrals in (\ref{eqn:MPMweak}) to be approximated by an assembly over the material points   
\begin{equation}\label{eqn:MPMweakDiscrete}
  \A_{ \forall p} \Bigl( [\nabla_x S_{vp}]^{T}\{\sigma_p\} V_p - [S_{vp}]^{T}\{b\} V_p + [S_{vp}]^{T}\{\dot{v}_p\} m_p \Bigr) - \A_{p\in P_c}\left(\{F_{N,vp} ^{\partial\Omega}\} + \{F_{T,vp} ^{\partial\Omega}\}\right) = \{0\}.
\end{equation}
Note that the subscripts on the basis functions, $[S_{vp}]$, have changed to highlight that they depend on both basis functions associated with the background grid and the characteristic function associated with the material point, $p$, which describes a material point's \emph{influence} (see Section~\ref{sec:basis} for details). The subscripts of the contact forces, $\{F_{N,vp} ^{\partial\Omega}\}$ and $\{F_{T,vp} ^{\partial\Omega}\}$, have also changed to highlight that contact is occurring at the material point $p\in P_c$, where $P_c$ is the set of material points in contact with the rigid body. 
In this paper it is assumed that Dirichlet boundary conditions are imposed directly on the nodes of the background grid as the imposition of general Dirichlet constraints is a separate area of research (see for example, \cite{Chandra2021,Bing2019,Singer2024,cortis2018imposition}). The traction term will be discretised in Section \ref{sec:contact}.

\subsection{Basis functions}\label{sec:basis}

In this paper we adopt the Generalised Interpolation Material Point (GIMP, \cite{bardenhagen2004generalized}) basis functions in order to mitigate the widely documented cell crossing instability and also to provide a convenient way to track contact between the deformable material (represented by material points) and a rigid body.  The basis functions, $S_{vp}$, are obtained by integrating the product of the shape functions associated with the background grid, $S_v$, with a characteristic function, $\chi_p$,  over the cuboid defining the domain associated with the material point, $\Omega_p$, normalised by the volume of the domain
\begin{equation}\label{eqn:Svp}
    S_{vp} = \frac{1}{V_d} \int_{\Omega_p} \chi_p S_v dV.  
\end{equation}
As elsewhere in the MPM literature, here $\chi_p$ is taken to be a unity hat function. The convolution in (\ref{eqn:Svp}) generates $C^1$ continuous basis functions that smooth the transfer of internal force from the material points to the grid nodes as a material point transitions between  elements.  It also, despite the name, generates basis functions that do not interpolate at the grid nodal locations unless $\chi_p$ is taken to be a Dirac delta function, and $S_{vp}=S_v$.  Note that the integral is normalised by the domain volume rather than the volume associated with the material point as, depending on the domain updating procedure, it is not guaranteed that the two volumes will be equal (see \citet{Coombs2020continuum} for a detailed discussion on domain updating procedures).  The domain updating procedure adopted in this paper is detailed in Section~\ref{sec:domain}.

\subsection{Discretisation in time}\label{eq:solution in time}

The discrete weak statement of equilibrium, (\ref{eqn:MPMweakDiscrete}), is discretised in time by a Newmark time integration scheme.  The nodal velocities are discretised as
\begin{equation}
    \{v^{n+1}_v\} = \frac{\gamma}{\Delta t\beta}\Bigl\{\{u_v^{n+1}\}-\{u_v^n\}\Bigr\}
    +\left(1-\frac{\gamma}{\beta}\right)\{v_v^n\}
    + \Delta t\left(\frac{1-\gamma}{2\beta}\right)\{\dot{v}_v^n\},
\end{equation}
where the superscripts $(\cdot)^{n+1}$ and $(\cdot)^{n}$ indicate quantities associated with the current and previous time step, respectively.  $\Delta t$ is the time increment and $\gamma$ and $\beta$ are the Newmark parameters.  The nodal accelerations are discretised as
\begin{equation}
    \{\dot{v}_v^{n+1}\} = \frac{1}{\Delta t^2\beta}\Bigl\{\{u_v^{n+1}\}-\{u_v^n\}\Bigr\}-\frac{1}{\beta\Delta t}\{v_v^n\}-\left(\frac{1}{2\beta}-1\right)\{\dot{v}_v^n\}.
\end{equation}
Note that for the MPM the nodal displacements at the end of each time step are discarded and the background grid reset or replaced.  Therefore it is assumed that $\{u^n_v\}=\{0\}$. The nodal velocities and accelerations at the start of the time step are projected from the material point values via
\begin{equation}
    \{v_v^n\} = [M]^{-1}\A_{ \forall p} \Bigl( [S_{vp}]^{T}\{{v}^n_p\} m_p \Bigr)
    \qquad \text{and} \qquad
    \{\dot{v}_v^n\} = [M]^{-1}\A_{ \forall p} \Bigl( [S_{vp}]^{T}\{\dot{v}^n_p\} m_p \Bigr),
\end{equation}
where $[M]=\A_{ \forall p} \left( [S_{vp}]^{T}[S_{vp}] m_p \right)$ is the consistent mass matrix. $\{{v}^n_p\}$ and $\{\dot{v}^n_p\}$ at the material point velocities and accelerations at the start of the time step (assumed to be equal to those from the end of the previous time step). The iterative increment in nodal displacements, $\{\delta u_v\}$, and the total displacement associated with a time step are
\begin{equation}
    \{\delta u_v\}  =  \left([K] +\frac{1}{\beta\Delta t^2}[M]\right)^{-1}\{f_{R}\}
    \qquad \text{and} \qquad
    \{u_v^{n+1}\} = \sum \{\delta u_v\},
\end{equation}
where $[K]=\A_{ \forall p} \left( [\nabla_x S_{vp}]^{T}[a][\nabla_x S_{vp}] V_p \right) + [K_c]$ is the consistent tangent stiffness matrix assembled at the nodes of the background grid, $[K_c]$ is the linearisation of the contact forces, $[a]$ is the spatial consistent tangent at the material point and $\{f_{R}\}$ is the out of balance force $R$esidual associated with (\ref{eqn:MPMweakDiscrete}), which should be equal to zero (within a tolerance) at the converged state. Iterative updates of the displacement, and other associated variables, continue until this condition is met. 
In this paper the Newmark parameters are taken to be $\gamma = 1/2$ and $\beta = 1/4$, resulting in an unconditionally stable implicit algorithm for linear problems \cite{Crisfield1997}.

\subsection{Stabilisation}


Since the work of \citet{Sulsky1995} it has been acknowledged that the MPM can suffer from instabilities linked to the arbitrary nature of the interaction between the material points and the background grid potentially leading to very small nodal masses in dynamic problems.  These small mass contributions result in conditioning issues and therefore issues with the inversion of the consistent mass matrix.  One way to avoid this problem, often done in explicit MPM implementations, is to lump the mass matrix.  However, this does not remedy the problem as the internal force calculation is linked to the derivatives of the basis functions, rather than the basis functions themselves, and therefore these forces on the mesh boundary are not necessarily small.  This can lead to very large spurious accelerations at the boundary of the body.  This problem is known as the ``small-cut'' issue in unfitted finite element methods \cite{Burman2010,Sticko2020}.  Some small mass cut off algorithms have been proposed \cite{Sulsky1995,Ma2010}, but these require specification of an arbitrary mass cut off threshold which is difficult to specify \cite{Sulsky1995}. The MUSL approach \cite{Sulsky1995} is a convenient way to mitigate the issue for explicit dynamic implementations but it cannot be applied to implicit approaches.  It also does not fix the conditioning issues of the mass matrix.  The stiffness matrix is also susceptible to ill conditioning issues when GIMP or B-spline basis functions are adopted due to their larger stencil and potential to generate small stiffness contributions to nodes near the boundary of the physical body \cite{coombs2022ghost}.  The extended B-spline approach of \citet{Yamaguchi2021} is a promising solution for B-spline basis functions, but it can't be applied to other basis functions.  
In this paper the ghost penalty approach of Burman \cite{Burman2010} is adopted to mitigate this issue based on the MPM implementation of Coombs \cite{coombs2022ghost}.  The method can be applied to explicit and implicit implementations of dynamic and quasi-static problems using any basis function.  The technique adds a stabilisation term into the mass and/or stiffness matrix that introduces additional continuity of the gradient of the solution across faces of the background grid at the boundary of the physical body (i.e. in the elements of the background grid that have the potential to be experiencing a small cut).

For linear elements the bi-linear form of the ghost stabilisation term\footnote{Note that the exponent on $h$ can be varied depending on the type of boundary condition at the edge of the physical domain, with an exponent of $1$ suggested for boundaries with Dirichlet constraints and $3$ for Neumann constraints \cite{Sticko2020}.  In this paper ghost stabilisation is primarily applied to free surfaces and surfaces subject to forces (i.e. tractions) from rigid body interactions.} is \cite{Sticko2020}
\begin{equation}\label{eqn:stabLinear}
	j(u_i,\eta_i) = \frac{h_f^{3}}{3} \int_{\Gamma} [[\partial_n u_i]]~[[\partial_n \eta_i]] d\Gamma,
    \qquad \text{where} \qquad
    [[u_i]] = u_i|_{F^+} - u_i|_{F^-}
\end{equation}
denotes the jump of $u_i$ across the face of two adjacent elements, with $F^+$ and $F^-$ being the faces of the \emph{positive} and \emph{negative} elements attached to the stabilised face. $\Gamma$ are the element faces where the stabilisation is applied, $n$ denotes the normal to the positive face and $h_f$ is a measure of the length of the face (usually taken to be the maximum side length).  The standard process of introducing the finite element approximation for the test and trial functions and eliminating the nodal values of the test function results in a stabilisation matrix with the form 
\begin{equation}\label{eqn:Jmatrix}
	[J_G] = \frac{h_f^{3}}{3}  \int_{\Gamma} \Bigl([G]^T\Bigl[[n][n]^T\Bigr][G]\Bigr) d\Gamma,
	\qquad \text{where} \qquad
	[G] = \Bigl[ [G^+] \quad -[G^-] \Bigr]
\end{equation}
contains the derivatives of the basis functions for the positive and negative elements and $[n]$ contains information on the normal direction of the face (see \cite{coombs2022ghost} for details and a full derivation).  The stabilisation contribution to the mass and stiffness matrices can be expressed as
\begin{equation}
    [M_G] = \gamma_M [J_G] \qquad \text{and} \qquad
    [K_G] = \frac{\gamma_K}{h^2} [J_G],
\end{equation}
where $\gamma_M$ and $\gamma_K$ are penalty parameters that control the degree of stabilisation added to the mass and stiffness matrices, respectively.  It is possible to link these parameters to the density and elastic stiffness of the material - see \citet{Sticko2020} and \citet{coombs2022ghost} for detailed discussions. In this paper $\gamma_M=\rho/4$ and $\gamma_K=E/30$, where $\rho$ and $E$ are the volume weighted average density and Young's modulus of the material points that occupy the elements that share the interface $\Gamma$.

\section{Rigid body}\label{sec: rigid body}
In order to describe the contact kinematics between the rigid body and the material points it is necessary to define various coordinate systems that will  as well as the allowable kinematics of the rigid body. Much of the details presented here can be found in the book by Wriggers \cite{wriggers2006computational}. However the mathematical tools which are used here are briefly detailed for both readability and a consistent nomenclature.

\subsection{Coordinate systems}\label{equ:rb coord system}
\begin{figure}[ht!]
    \centering
    \def\svgwidth{\textwidth}
   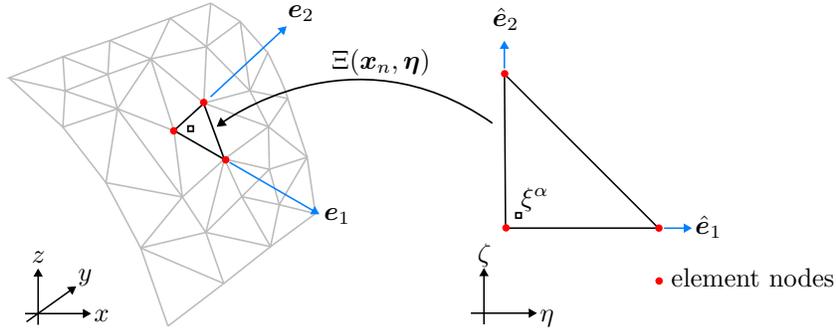
    \caption{The rigid body global and local coordinate systems.}
    \label{fig:rigid body coordiantes}
\end{figure}

The surface of the rigid body surface is discretised by a triangular mesh, as in the left of Figure \ref{fig:rigid body coordiantes}.
On each triangle in the global domain, ${K}$, a contravariant coordinate system is defined, $\bm{x}^\prime = \bm{e}_\alpha\xi^\alpha$, where $\bm{x}^\prime\in K$, is shown by the square marker in Figure \ref{fig:rigid body coordiantes}. Similarly on the local domain $\widehat{K}$, shown on the right Figure \ref{fig:rigid body coordiantes}, another contravariant coordinate system is defined, $\bm{\eta} = \hat{\bm{e}}_\alpha\xi^\alpha$. Between the two coordinate systems there is the affine mapping, 
\begin{equation}\label{equ: aff mapping}
    \bm{x}^\prime  = \bm{\Xi}(\bm{x}_n,\bm{\eta})
\end{equation}
where $\bm{x}_n$ are the nodal positions of the rigid body triangle. The mapping is such that the contravariant components, $\xi^\alpha$, are convected between the coordinate systems \cite{dvorkin2006nonlinear},
\begin{equation}\label{eq:x map}
    \bm{x}^\prime={\bm{e}}_\alpha\xi^\alpha 
    \quad\text{and}\quad
    \bm{\eta}=\hat{\bm{e}}_\alpha\xi^\alpha
\end{equation}
such that the value $\xi^\alpha$ remains constant between the local and global domains when the mapping \eqref{equ: aff mapping} is performed. In order to track the movement of a contact point over a triangle, it is convenient to do so with nodal shape functions and the local coordinates $\xi^\alpha$ defined at the nodes $\xi^\alpha_n$, where $n$ is the nodal number. The value of $\xi^\alpha$ can therefore be described as
\begin{equation}\label{eq:xi sf}
    \xi^\alpha = \sum_n^3N_n(\xi)\xi^\alpha_n
\end{equation}
where $N_n(\xi)$ are the nodal shape functions. Substituting Equation \eqref{eq:xi sf} into Equation \eqref{eq:x map} gives the description for the position $\bm{x}^\prime$ but using the local element shape functions,
 \begin{equation}\label{eq:tri shape functions}
      \bm{x}^\prime=\bm{e}_\alpha\sum_n^3N_n(\xi)\xi^\alpha_n=\sum_n^3N_n(\xi)\bm{x}_n
 \end{equation}
where $\bm{x}_n$ is the position of the triangle corners in the global domain. Additionally, the tangent to the rigid body in the global domain $\bm{t}_\alpha$, see \citet{pietrzak1999large}, is found with the derivative of shape functions with respect to $\xi^\alpha$
 \begin{equation}
      \bm{t}_\alpha = \frac{\partial\bm{x}}{\partial\xi^\alpha}=\bm{e}_\beta\sum_n^3\frac{\partial N_n(\xi)}{\partial \xi^\alpha}\xi^\beta_n
                    = \sum_n^3\frac{\partial N_n(\xi)}{\partial \xi^\alpha}\bm{x}_n
 \end{equation}

\subsection{Kinematics}

\begin{figure}[ht!]
    \centering
     \def\svgwidth{0.6\textwidth}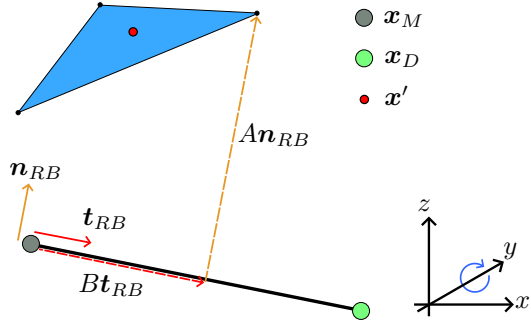
        \caption{Rigid body: truss frame description.}\label{fig: RB truss main}
\end{figure}

A rigid body is defined by a set of notes unable to move relative to each other, the motion of all the nodes is described with a single set of rigid body degrees of freedom. 
The work presented here was originally developed for geotechnical engineering problems which often involve a chain or wire pulling a rigid body, such as a plough or anchor. Here a truss frame is used to model the chain/wire, and the frame is also used to describe the kinematics of the rigid body. The kinematic description is achieved by using the vertices of a truss element to describe the position of all nodes. An example of this for a single triangle of the rigid body is shown in Figure \ref{fig: RB truss main}, where the motion of the truss frame is restricted to the $y$-plane.
The truss element in Figure \ref{fig: RB truss main} is comprised of two nodes, $\bm{x}_{M}$ and $\bm{x}_{D}$ from which a tangent $\bm{t}_{rb}$ and normal $\bm{n}_{RB}$ are defined. They are respectively computed by 
\begin{equation}
    \bm{t}_{RB} = \frac{\bm{x}_{D} - \bm{x}_{M}}{\|\bm{x}_{D} - \bm{x}_{M}\|}\quad\text{and}\quad \bm{n}_{RB} = \bm{R}\cdot\bm{t}_{RB}
\end{equation}
where $\bm{R}$ is a rotation matrix of $90^\circ$ about the $y$-axis. The nodal position of the triangle element, $\bm{x}_n$, can now be described as function of $\bm{x}_{M}$, $\bm{t}_{RB}$, $\bm{n}_{RB}$ and two scalar values $B_n$ and $A_n$ which have unique, and constant, values for each node,
\begin{equation}\label{equ:app xn definition m}
    \bm{x}_n = A_n\bm{I}\cdot\bm{n}_{RB} + B_n\bm{I}\cdot\bm{t}_{RB} + \bm{x}_{M} = (A_n\bm{R} + B_n\bm{I})\cdot\bm{t}_{RB} + \bm{x}_{M} = \bm{G}_n(\bm{\theta})
\end{equation}
where $\bm{I}$ is the identity matrix and $\bm{G}_n(\bm{\theta})$ is the function for the position $\bm{x}_n$ dependent on the variable $\bm{\theta}$ which is a vector containing the points $\bm{x}_{M}$ and $\bm{x}_{D}$. Every node in the rigid body is described by Equation \eqref{equ:app xn definition m} which for all nodes has the variables $\bm{x}_{M}$ and $\bm{x}_{D}$, hence there is potential for relative motion between the nodes. Substituting Equation \eqref{equ:app xn definition m} into \eqref{eq:tri shape functions} gives,
\begin{equation}\label{equ:app Nx full m}
    \bm{x}^\prime = \sum_n^3N_n(\xi)\left[(A_n\bm{R} + B_n\bm{I}):\bm{t}_{RB} + \bm{x}_{M}\right]= \sum_n^3N_n(\xi)\bm{G}_n(\bm{\theta}).
\end{equation}
In Section \ref{sec:contact} $\bm{x}^\prime$ is used as the contact point of the rigid body with the material point, for which the implicit solver will require its first and second variations, provided in Appendix \ref{App:truss frame}.

Lastly, the residual for the equation of linear momentum of the truss frame, with the contact forces, is including in the monolithic solver, presented in Section \ref{eq:solution in time}, to form a coupled problem. The residual for linear momentum takes the form,
\begin{equation}\label{eqn:frameweak}
    \A_{f\in F}\{\sigma_f\} \text{d}V 
  - \A_{f\in F}\{b\} \text{d}V 
  + \A_{f\in F}[m_f]\{\dot{v_f}\} \text{d}V
  - \A_{l\in L}\left(\{F_N^{l}\} + \{F_T^{l}\}\right)\text{d}s
  =\{0\},
\end{equation}
where $F$ is the set for all truss elements $f$, $L$ is the set for all rigid body triangles $l$, $[m]$ is the nodal mass matrix for $f$,  $\{b_f\}$ is the nodal external force vector, $\{\dot{v_f}\}$ is the nodal acceleration and the last term is the contact forces acting on the triangles of the rigid body, split into normal, $\{F_N^{l}\}$ and tangential components $\{F_T^{l}\}$. $\{\sigma_f\}$ is calculated as
\begin{equation*}
\begin{split}
\{\sigma_f\}  = F\left[ -c,  -s,  c,  s \right]^T&,\quad\text{where}\quad
F  = E_f (L_0 - L),\quad
L  = \|\bm{x}_1 - \bm{x}_2\|,\\
c  &= \frac{x_2 - x_1}{L}\quad\text{and}\quad
s  = \frac{z_2 - z_1}{L}
\end{split}
\end{equation*}
where $L_0$ is the original length and the subscripts $1$ and $2$ correspond to the two nodes of the truss element.

\section{Contact}\label{sec:contact}
This section presents the description of the contact forces present in Equations \eqref{eqn:MPMweakDiscrete} and \eqref{eqn:frameweak}. When considering the contact between the material point $p$ and the triangle $l$ the general form for the contact forces is 
\begin{equation}
\delta U = \delta U_N + \delta U_T = 
\int_{\psi_t(l)\cap{\delta}(\bm{x}^\prime)}
\left(\{F_N^{l}\} + \{F_{N,vp} ^{\partial\Omega}\}\right)\text{d}s
+  \int_{\psi_t(l)\cap{\delta}(\bm{x}^\prime)}
\left(\{F_T^{l}\} + \{F_{T,vp} ^{\partial\Omega}\}\right)\text{d}s
\end{equation}
where $\delta(\bm{x}^\prime)$ is the Dirac delta function and $\bm{x}^\prime$ is the projection of the position of the material point $p$ onto the rigid body surface, using the Closest Point Projection (CPP) scheme, detailed in Section \ref{sec:gap function}.

\subsection{Contact introduction}
This section describes the contact formulation between the rigid body and the GIMPs. Using the same contact methods as in FEM, the contact formulation and methodology are driven by the gap function \cite{wriggers2006computational,curnier1995continuum}. The gap function is the projection of a point from the secondary surface (s), with coordinates $\bm{x}^\prime$, onto the main surface (m), with coordinates $\bm{x}$. The gap function is used to: i) determine if there is contact and between which point and surface; ii) how much overlap there is; and iii) through the rate derivative of the gap function to determine the relative tangential velocity.

Here the CPP scheme is used to determine the gap function, see the works of Curnier and their coworkers \cite{curnier1995continuum,pietrzak1999large}. However, other methods do exists such as ray-tracing (see \citet{poulios2015unconstrained}), which has been applied to GIMP-to-GIMP contact, \citet{pretti2024continuum}. 

\subsection{Gap function}\label{sec:gap function}
In order to define the gap function, first the definition of the main (m) and secondary (s) surfaces/domains need to be established, as well as the general contact methodology. Here a point-to-surface contact is used, where the points are the primary domain, consisting of the GIMP vertices, and the secondary domain is the triangular surface of the rigid body.

The CPP gap function is the minimum distance of a point's projection point onto the surface, defined as,
\begin{equation}\label{equ:gap function}
    \bm{g}_N(\xi(\tau)^\alpha,\tau) = g_N(\xi(\tau)^\alpha,\tau)\bm{n}(\xi(\tau)^\alpha,\tau) = 
    \bm{x}(\tau) - \bm{x}^\prime(\xi(\tau)^\alpha,\tau)\qquad\text{and}\qquad g_N = \bm{g}_N\cdot \bm{n},
\end{equation}
where $\tau$ is time and $\xi(\tau)^\alpha$ is a function that describes that the local location of the projection on the surface. 
The normal contact law is described by the Signorini-Hertz-Moreau conditions, which are a function of $g_N$ and the corresponding normal penetration force $p_N$,
\begin{equation}
    g_N \left\{
    \begin{tabular}{cl}
        $=0$ & contact \\
        $\geq 0$ & no contact
    \end{tabular}
    \right.
    ,\quad
    p_N \left\{
    \begin{tabular}{cl}
        $\leq 0$ & contact \\
        $= 0$ & no contact
    \end{tabular}
    \right.
    ,\quad
    \text{and}
    \quad
    g_Np_N = 0.
\end{equation}
Where the last condition is relaxed for the penalty method.

Following the work of \citet{curnier1995continuum} an objective definition of the relative velocities can be obtained as a total time derivative of $\bm{g}_N$, which is necessary in order to evaluate the frictional component of contact. It is expressed as 
\begin{equation}\label{equ: gap time div contact}
\begin{split}
     \dot{\bm{g}}_N(\xi(\tau)^\alpha,\tau) 
    &= 
    \dot{g}_N(\xi(\tau)^\alpha,\tau)\bm{n}(\xi(\tau)^\alpha,\tau)
    +
    g_N(\xi(\tau)^\alpha,\tau) \dot{\bm{n}}(\xi(\tau)^\alpha,\tau)\\
    :&=\bm{\mathring{g}_N}(\xi(\tau)^\alpha,\tau) + \bm{\mathring{g}_t}(\xi(\tau)^\alpha,\tau)\\  
    :&= \bm{\mathring{g}_N}(\xi(\tau)^\alpha,\tau) +  \bm{t}_\alpha(\xi(\tau)^\alpha,\tau)\dot{\xi}(\tau)^\alpha
\end{split}
\end{equation}
where $\bm{\mathring{g}_N}$ is rate of penetration and $\bm{\mathring{g}_t}$ is the relative tangential velocity between the two points. Following the work of \cite{wriggers2006computational,curnier1995continuum} the residual for contact can be expressed as 
\begin{equation}\label{eq:d contact energy}
    \delta U = \int_{\psi_t(l)\cap{\delta}(\bm{x}^\prime)}
     \left(\delta{g}_N\bm{n}\cdot\bm{p}_N
    + 
    \delta{\xi}^\alpha\bm{t}_\alpha\cdot\bm{p}_t \right)\text{ d}x
\end{equation}
where ${\delta}(\bm{x}^\prime)$ is the Dirac delta function at the current position of the secondary surface. Additionally the first variation of the normal and tangential gap rates are required
\begin{equation}\label{equ: gap variations}
    \bm{\mathring{g}_N} \delta \tau = \bm{n}\dot{g}_N\delta \tau = \bm{n}\delta{g}_N
    \qquad \text{and} \qquad 
    \bm{\mathring{g}_t} \delta \tau = \bm{t}_\alpha\dot{\xi}^\alpha\delta \tau = \bm{t}_\alpha\delta{\xi}^\alpha.
\end{equation}

\subsection{Normal contact}
The normal component of the residual, Equation \eqref{eq:d contact energy}, is
\begin{equation}\label{equ: d contact energy normal}
      \delta U_N = \int_{\psi_t(l)\cap{\delta}(\bm{x}^\prime)}
     \left(\delta{g}_N\bm{n}\cdot\bm{p}_N \right)\text{ d}x 
\end{equation}
where $\bm{p}_N$ is the normal pressure that resists contact. Here the penalty method is used to calculate $\bm{p}_N$, which is defined in terms of a normal penalty parameter $\epsilon_N$ and the gap function,
\begin{equation}\label{equ: penalty function}
    \bm{p}_N = \epsilon_N\bm{n}g_N
\end{equation}
such that \eqref{equ: d contact energy normal} becomes
\begin{equation}\label{equ: d contact energy normal expanded}
      \delta U_N = \int_{\psi_t(l)\cap{\delta}(\bm{x}^\prime)}
     \left(\delta{g}_N\bm{n}\cdot\epsilon_N\bm{n}g_N \right)\text{ d}x 
\end{equation}
%
Linearising \eqref{equ: d contact energy normal expanded} gives the gradient of the residual, necessary for the Newton-Raphson scheme,
\begin{equation}\label{equ: d contact energy normal K}
      \Delta\delta U_N = \int_{\psi_t(l)\cap{\delta}(\bm{x}^\prime)}
     \left(\Delta(\delta{g}_N\bm{n})\cdot\epsilon_N\bm{n}g_N
     +
     \delta{g}_N\bm{n}\cdot\Delta(\epsilon_N\bm{n}g_N) \right)\text{ d}x,
\end{equation}
where $\Delta(\delta{g}_N\bm{n})$ and $\Delta(\delta{g}_N\bm{n})$ are defined in Appendix \ref{App:gap var}.

\subsection{Tangential contact}
The frictional component of the residual, Equation \eqref{eq:d contact energy}, is
\begin{equation}\label{equ: d contact energy tangent}
      \delta U_T = \int_{\psi_t(l)\cap{\delta}(\bm{x}^\prime)}
     \left(\delta{\xi}^\alpha\bm{t}_\alpha\cdot\bm{p}_T  \right)\text{ d}x 
\end{equation}
where $\bm{p}_T$ is the tangential pressure that acts to resist sliding motion contact. Similarly to the normal contact, linearising \eqref{equ: d contact energy tangent} gives 
\begin{equation}\label{equ: d contact energy tangent K}
      \Delta\delta U_T = \int_{\psi_t(l)\cap{\delta}(\bm{x}^\prime)}
     \left(\Delta(\delta{\xi}^\alpha)\bm{t}_\alpha\cdot\bm{p}_T
     +
     \delta{\xi}^\alpha\Delta(\bm{t}_\alpha\cdot\bm{p}_T)) \right)\text{ d}x 
\end{equation}
where $\Delta(\delta{\xi}^\alpha)$ and $\Delta(\bm{t}_\alpha\cdot\bm{p}_T)$ are defined in Appendix \ref{App:xi var}.

\subsection{Friction}\label{sec: friction}
Here the methodology for determining the friction acting between the two bodies is presented. It is an extension of the 2D method proposed by \citet{bird2024implicit} for GIMP-to-rigid body contact, which itself is based on the work of \citet{wriggers2006computational}. The methodology starts by defining the total tangential movement of a point on the surface, 
\begin{equation}
    \bm{g}_T = \int_\tau \mathring{\bm{g}_T}\text{ d}t= \int_\tau \mathring{\bm{g}}_{slip}\text{ d}t +  \int_\tau \mathring{\bm{g}}_{stick}\text{ d}t = \bm{g}_{slip} + \bm{g}_{stick},
\end{equation}
which is comprised of the components, $\bm{g}_{slip}$, which is the purely dissipative and forms the plasticity component of the friction model, and (ii) the tangential stick, $\bm{g}_{stick}$. 

Here the Coulomb friction law is used where the frictional force is defined as,
\begin{equation}\label{eq: slip force}
   \bm{p}_T = \mu|p_N|\frac{\mathring{\bm{g}}_{slip}}{||\mathring{\bm{g}}_{slip}||}\qquad\text{if}\qquad ||\bm{p}_{stick}||>\mu|p_N|,
\end{equation}
where $\mu$ is the coefficient of friction, which is assumed constant. The friction law is subject to the Karush-Kuhn-Tucker (KKT) conditions,
\begin{equation}\label{eq:KKT condition}
    f =||\bm{p}_{stick}||-\mu p_N \leq 0,\qquad \lambda\geq 0 \qquad \text{and} \qquad f\lambda=0,
\end{equation}
where $f$ is the frictional yield function, $\lambda$ is the yield rate. Here we use an elastic law, or penalty method, to describe the force corresponding to $\bm{g}_{stick}$, therefore $\bm{g}_{stick}$ can be described as purely elastic and a recoverable motion of the particle along the surface
\begin{equation}
    \bm{p}_{stick} = \epsilon_T\bm{g}_{stick}.
\end{equation}

To implement the continuous friction laws within the material point framework, the equations need to be reformulated into the quasi-static problem and solved using the same methods used to resolve elastic-perfectly-plastic material laws. First, the total tangential movement over the time step $m+1$ is defined,
\begin{equation}\label{eq tangental slip}
    \begin{split}
        \mathring{\bm{g}_T}\Delta t &\approx {\bm{g}_{T,m+1} - \bm{g}^{m}_{T,m}} \\
        &= \bm{t}_{\alpha,m+1} \Delta{\xi^\alpha}= \bm{t}_{\alpha,m+1}({\xi}_{m+1}^\alpha-{\xi}_{m}^\alpha) \\
        &= \Delta \bm{g}_T. 
    \end{split}
\end{equation} 
Next the elastic trial stick state can be defined as 
\begin{equation}
\bm{p}_{tr}=\epsilon_T(\bm{g}_{T,m+1} -\bm{g}_{slip,m})
\end{equation}
where $(\bm{g}_{T,m+1} -\bm{g}_{slip,m})$ is the trial stick movement that has occurred. However, since the friction coefficient $\mu$ in Equation \eqref{eq:KKT condition} is constant, and considering the arguments presented in \cite{bird2024implicit,wriggers2006computational}, the elastic trial state can be expressed as
\begin{equation}
\bm{p}_{tr}=\bm{p}_{T,m} + \epsilon_T\Delta \bm{g}_T,
\end{equation}
where $\bm{p}_{T,m}$ is the previously converged tangential friction vector. With these definitions, the trial stick state is used with yield surface \eqref{eq:KKT condition} to determine if the contact is in a stick or slip state, and subsequently the frictional force can be calculated,
\begin{equation}\label{equ:stick slip}
    \bm{p}_T = \left\{ \begin{tabular}{cccl}
          $\bm{p}_{tr}$ & if & $f\leq 0$, & stick \\
          $\mu|p_N|\left({{\bm{p}_{tr}}}/{||{\bm{p}_{tr}}||}\right) $ & if & $f> 0$, & slip
    \end{tabular}\right.
\end{equation}
The variation of $\bm{p}_T$ is required for the monolithic solver, as shown in Appendix \ref{App:stick and slip}.

\section{Numerical implementation}
This section details two critical elements of the material point-rigid body interaction that need careful treatment in order to realise and accurate and robust numerical implementation, namely: (i) the contact detection approach using the corners of the GIMP domains and (ii) the update of these domains to maintain a stable, energy consistent solution.  

\subsection{Closest Point Projection (CPP) for contact detection}\label{sec: CPP}
The description of the gap function is fundamental to modelling contact, as described in Section \ref{sec:contact}. However, when considering contact between a point and an irregular surface there are algorithmic nuances that need to be discussed which can be overlooked when going from a mathematical-to-algorithmic description, particular in the case with GIMPs. The focus of this section is the calculation of Equation \eqref{equ:gap function} and presents an extension of the 2D contact methodology presented by \citet{bird2024implicit} to 3D.

\begin{figure}[ht!]
    \centering
     \def\svgwidth{0.4\textwidth}
\begingroup%
  \makeatletter%
  \providecommand\color[2][]{%
    \errmessage{(Inkscape) Color is used for the text in Inkscape, but the package 'color.sty' is not loaded}%
    \renewcommand\color[2][]{}%
  }%
  \providecommand\transparent[1]{%
    \errmessage{(Inkscape) Transparency is used (non-zero) for the text in Inkscape, but the package 'transparent.sty' is not loaded}%
    \renewcommand\transparent[1]{}%
  }%
  \providecommand\rotatebox[2]{#2}%
  \newcommand*\fsize{\dimexpr\f@size pt\relax}%
  \newcommand*\lineheight[1]{\fontsize{\fsize}{#1\fsize}\selectfont}%
  \ifx\svgwidth\undefined%
    \setlength{\unitlength}{252.16014952bp}%
    \ifx\svgscale\undefined%
      \relax%
    \else%
      \setlength{\unitlength}{\unitlength * \real{\svgscale}}%
    \fi%
  \else%
    \setlength{\unitlength}{\svgwidth}%
  \fi%
  \global\let\svgwidth\undefined%
  \global\let\svgscale\undefined%
  \makeatother%
  \begin{picture}(1,1.3212816)%
    \lineheight{1}%
    \setlength\tabcolsep{0pt}%
    \put(0,0){\includegraphics[width=\unitlength,page=1]{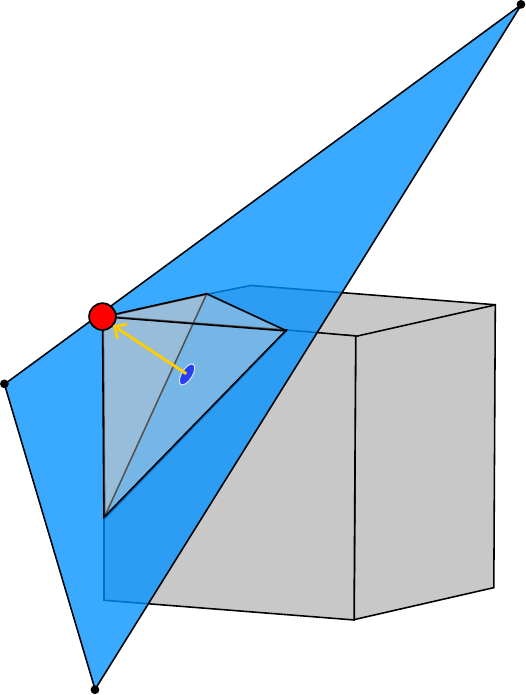}}%
    \put(0.13194932,0.76609296){\makebox(0,0)[lt]{\lineheight{1.25}\smash{\begin{tabular}[t]{l}$\bm{x}_{gn}$\end{tabular}}}}%
    \put(0,0){\includegraphics[width=\unitlength,page=2]{Contact_triangle.pdf}}%
    \put(0.37306505,0.59402259){\makebox(0,0)[lt]{\lineheight{1.25}\smash{\begin{tabular}[t]{l}$\bm{x}^\prime$\end{tabular}}}}%
    \put(0,0){\includegraphics[width=\unitlength,page=3]{Contact_triangle.pdf}}%
    \put(0.24166547,0.61234723){\makebox(0,0)[lt]{\lineheight{1.25}\smash{\begin{tabular}[t]{l}$\bm{g}$\end{tabular}}}}%
  \end{picture}%
\endgroup%

        \caption{CPP: A contact vertex (red dot) of a GIMP domain, shown in grey, in contact with a rigid body triangle, blue. The gap function is also shown by the yellow arrow.}\label{fig: GIMP contact p-to-tri}
    \label{fig:enter-label}
\end{figure}

The vertices of the GIMP domains are the points in contact with the faceted surface of the rigid body, consistent with the point-to-surface contact described in Section \ref{sec:contact}. This is necessary so that contact occurs on the material boundary, otherwise the contact is not consistent and spurious stresses are observed for the contact GIMPs \cite{bird2024implicit}. Since a GIMP defines the material domain over which it is integrated there is an explicit representation of the boundary, this can be observed through the definition of the GIMP basis (see Equation (\ref{eqn:Svp})). The integration of (\ref{eqn:Svp}) to construct ${S}_{vp}$ is exact, thus through inspection of calculation of the linear momentum residual, \eqref{eqn:MPMweakDiscrete}, is also an integral of a function over an exact domain. The point to surface contact between a vertex of the GIMP and a triangle of the rigid body surface is described with Figure \ref{fig: GIMP contact p-to-tri}
where, the position of the GIMP vertex is described as 
\begin{equation}
    \bm{x}_{gn} = \sum_{v \in E}N_v\bm{x}_v
\end{equation}
which when in contact is the point $\bm{x}$ that is projected onto point $\bm{x}^\prime$ in Equation \eqref{equ:gap function}. 

As described in Section \ref{sec: rigid body} the surface of the rigid body is constructed from a triangular mesh to produce a constant normal direction over each triangle surface. Points and lines of the mesh do not have a normal description and hence it is not possible to define a gap function for them, see the work of \citet{curnier1995continuum} for this argument. Hence contact is only observed on the triangle faces.

\begin{figure}[ht!]
    \centering
    \begin{subfigure}[b]{0.55\textwidth}
        \centering
         \includegraphics[width=\linewidth]{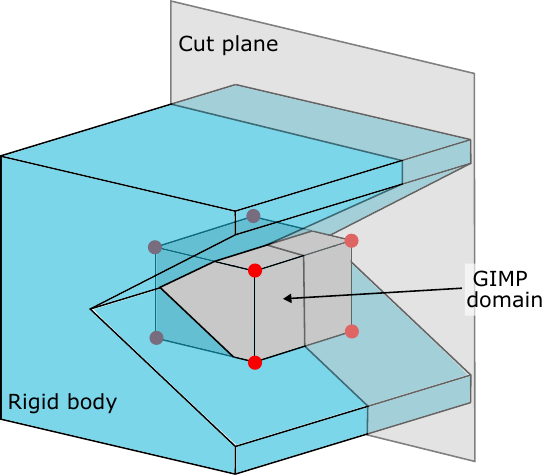}
        \caption{}\label{fig: contact regions 1}
    \end{subfigure}%
    \hfill
     \begin{subfigure}[b]{0.33\textwidth}
        \centering
         \def\svgwidth{\textwidth}
\begingroup%
  \makeatletter%
  \providecommand\color[2][]{%
    \errmessage{(Inkscape) Color is used for the text in Inkscape, but the package 'color.sty' is not loaded}%
    \renewcommand\color[2][]{}%
  }%
  \providecommand\transparent[1]{%
    \errmessage{(Inkscape) Transparency is used (non-zero) for the text in Inkscape, but the package 'transparent.sty' is not loaded}%
    \renewcommand\transparent[1]{}%
  }%
  \providecommand\rotatebox[2]{#2}%
  \newcommand*\fsize{\dimexpr\f@size pt\relax}%
  \newcommand*\lineheight[1]{\fontsize{\fsize}{#1\fsize}\selectfont}%
  \ifx\svgwidth\undefined%
    \setlength{\unitlength}{165.87608037bp}%
    \ifx\svgscale\undefined%
      \relax%
    \else%
      \setlength{\unitlength}{\unitlength * \real{\svgscale}}%
    \fi%
  \else%
    \setlength{\unitlength}{\svgwidth}%
  \fi%
  \global\let\svgwidth\undefined%
  \global\let\svgscale\undefined%
  \makeatother%
  \begin{picture}(1,0.85510298)%
    \lineheight{1}%
    \setlength\tabcolsep{0pt}%
    \put(0,0){\includegraphics[width=\unitlength,page=1]{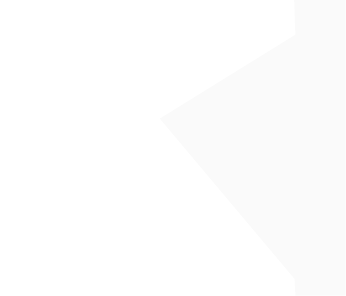}}%
    \put(0.93677764,0.69225478){\rotatebox{-90}{\makebox(0,0)[lt]{\lineheight{1.25}\smash{\begin{tabular}[t]{l}Outside rigid body\end{tabular}}}}}%
    \put(0,0){\includegraphics[width=\unitlength,page=2]{contact_2D_GIMP_contact.pdf}}%
    \put(0.00958104,0.01818569){\makebox(0,0)[lt]{\lineheight{1.25}\smash{\begin{tabular}[t]{l}Inside rigid body\end{tabular}}}}%
    \put(0,0){\includegraphics[width=\unitlength,page=3]{contact_2D_GIMP_contact.pdf}}%
    \put(0.38205894,0.4091011){\makebox(0,0)[lt]{\lineheight{1.25}\smash{\begin{tabular}[t]{l}$\bm{g}$\end{tabular}}}}%
    \put(0.5421299,0.47750254){\makebox(0,0)[lt]{\lineheight{1.25}\smash{\begin{tabular}[t]{l}$\bm{x}^\prime$\end{tabular}}}}%
    \put(0.26416638,0.19739188){\makebox(0,0)[lt]{\lineheight{1.25}\smash{\begin{tabular}[t]{l}$\bm{x}_{gn}$\end{tabular}}}}%
    \put(0,0){\includegraphics[width=\unitlength,page=4]{contact_2D_GIMP_contact.pdf}}%
  \end{picture}%
\endgroup%

        \caption{}\label{fig: contact regions 2D}
    \end{subfigure}%
    \caption{Closest point projection: (a) is the 3D view of a GIMP in contact with a concave edge of the rigid body, (b) is the corresponding slice view through the domain with green regions showing where the CPP is onto triangles, and red regions the CPP onto points/lines of the rigid body.}\label{fig: gap projections}
\end{figure}

An issue with not being able to define a gap function for the points/lines is that it may be possible for points to travel into the domain undetected. However this problem is mitigated by each GIMP having 8 potential points of contact, furthermore and points/lines can be filleted, in CAD software, to further prevent this phenomenon. A 3D scenario of where contact is not detected at some GIMP corners is shown in Figure \ref{fig: contact regions 1}, with a 2D slice shown in Figure \ref{fig: contact regions 2D}. The green and red region is where respectively the contact constraint can, and cannot, be applied. 
~
~
From an algorithmic point of the view these red/green zones are not identified explicitly; the red regions are simply the result of the CPP being onto a line or point, and the green regions a CPP onto a triangle.

\subsection{GIMP domain update procedure}\label{sec:domain}
\begin{figure}[ht!]
    \centering
    \begin{subfigure}[b]{0.2\textwidth}
        \centering
         \def\svgwidth{\textwidth}
\begingroup%
  \makeatletter%
  \providecommand\color[2][]{%
    \errmessage{(Inkscape) Color is used for the text in Inkscape, but the package 'color.sty' is not loaded}%
    \renewcommand\color[2][]{}%
  }%
  \providecommand\transparent[1]{%
    \errmessage{(Inkscape) Transparency is used (non-zero) for the text in Inkscape, but the package 'transparent.sty' is not loaded}%
    \renewcommand\transparent[1]{}%
  }%
  \providecommand\rotatebox[2]{#2}%
  \newcommand*\fsize{\dimexpr\f@size pt\relax}%
  \newcommand*\lineheight[1]{\fontsize{\fsize}{#1\fsize}\selectfont}%
  \ifx\svgwidth\undefined%
    \setlength{\unitlength}{70.86614173bp}%
    \ifx\svgscale\undefined%
      \relax%
    \else%
      \setlength{\unitlength}{\unitlength * \real{\svgscale}}%
    \fi%
  \else%
    \setlength{\unitlength}{\svgwidth}%
  \fi%
  \global\let\svgwidth\undefined%
  \global\let\svgscale\undefined%
  \makeatother%
  \begin{picture}(1,1.61018234)%
    \lineheight{1}%
    \setlength\tabcolsep{0pt}%
    \put(0,0){\includegraphics[width=\unitlength,page=1]{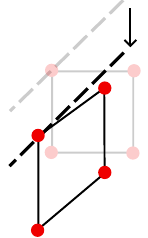}}%
    \put(0.92523562,1.44851122){\color[rgb]{0,0,0}\makebox(0,0)[lt]{\lineheight{1.25}\smash{\begin{tabular}[t]{l}$u^{rb}_i$\end{tabular}}}}%
  \end{picture}%
\endgroup%

        \caption{}\label{fig: GIMP contact}
    \end{subfigure}%
    \hspace{1cm}
    \begin{subfigure}[b]{0.2\textwidth}
        \centering
            \def\svgwidth{\textwidth}
\begingroup%
  \makeatletter%
  \providecommand\color[2][]{%
    \errmessage{(Inkscape) Color is used for the text in Inkscape, but the package 'color.sty' is not loaded}%
    \renewcommand\color[2][]{}%
  }%
  \providecommand\transparent[1]{%
    \errmessage{(Inkscape) Transparency is used (non-zero) for the text in Inkscape, but the package 'transparent.sty' is not loaded}%
    \renewcommand\transparent[1]{}%
  }%
  \providecommand\rotatebox[2]{#2}%
  \newcommand*\fsize{\dimexpr\f@size pt\relax}%
  \newcommand*\lineheight[1]{\fontsize{\fsize}{#1\fsize}\selectfont}%
  \ifx\svgwidth\undefined%
    \setlength{\unitlength}{70.86614173bp}%
    \ifx\svgscale\undefined%
      \relax%
    \else%
      \setlength{\unitlength}{\unitlength * \real{\svgscale}}%
    \fi%
  \else%
    \setlength{\unitlength}{\svgwidth}%
  \fi%
  \global\let\svgwidth\undefined%
  \global\let\svgscale\undefined%
  \makeatother%
  \begin{picture}(1,1.61018234)%
    \lineheight{1}%
    \setlength\tabcolsep{0pt}%
    \put(0,0){\includegraphics[width=\unitlength,page=1]{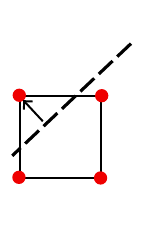}}%
    \put(0.24467669,0.87072534){\color[rgb]{0,0,0}\makebox(0,0)[lt]{\lineheight{1.25}\smash{\begin{tabular}[t]{l}$g_i$\end{tabular}}}}%
  \end{picture}%
\endgroup%

        \caption{}\label{fig: GIMP contact reset bad}
    \end{subfigure}
      \hspace{1cm} 
    \begin{subfigure}[b]{0.2\textwidth}
        \centering
            \def\svgwidth{\textwidth}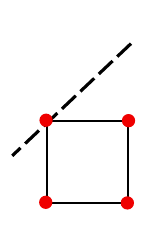
        \caption{}\label{fig: GIMP contact reset good}
    \end{subfigure}
    \caption{GIMP contact update: The converged position of the GIMP domain for time step $n$ is shown in (a), the updated position for the beginning of step $n+1$ is shown in (b) and the updated position  for the beginning of step $n+1$ after the gap minimisation in (c). The surface of the rigid body is shown by the dashed line}\label{fig:GIMP contact update intro}
\end{figure}

Large time steps can be taken when solving a contact problem implicitly, and in the GIMPs which are in contact can undergo significant deformation, as shown in Figure \ref{fig: GIMP contact}, during which the normal contact condition remains met. Once the time step is complete a reset algorithm is performed on the mesh and GIMP domains, see Figure \ref{fig: method steps f}, during which the GIMPs are reset from a distorted to regular shape, with the volume of the GIMP domain equal to the material volume, see Figure \ref{fig: GIMP contact reset bad}. Whilst this is acceptable in other regions of the mesh, for GIMPs in contact a large, purely numerical, increase in the gap function is introduced, as shown in Figure \ref{fig: GIMP contact reset bad}, despite no external work being performed. The result is that at the beginning of the new time step a significant unrealistic normal contact force can exist. This affects the numerical stability of the algorithm since a large non-linearity is introduced which often results in non-convergence of the  solver. This overlap is particularly problematic since the size of the non-linearity is the same, regardless of the size of the time step. Therefore reducing the time step size, and corresponding load increment size, can restore stability, but the solver will continue to fail regardless of the time step size, consequently the simulation will fail.

To solve this problem an energy minimisation is performed during the update of the GIMP. The energy minimisation is presented as, 
\begin{equation}\label{equ: gap energy minimisation}
\text{argmin}_{x_i \in \mathcal{R}^3} \, E(x_i)\quad\text{where}\quad  E(x_i) = \sum_{n}^8\left[\epsilon_{N,n}^{old}g_{N,n}^{old} - \epsilon_{N,n}^{new}g_{N,n}^{new}(x_i)\right]^2
\end{equation}
where $x_i$ is the position of the GIMP domain centre, $g_{N,n}$ is the magnitude of the normal gap function at node $n$, $\epsilon_{N,n}$ is the normal penalty, with the superscripts $old$ and $new$ corresponding to the previously converged and current variable state. The authors emphasise that this \emph{is} a solution to the problem, but there could be other possible variations on this perhaps also considering how much the material is moved to satisfy \eqref{equ: gap energy minimisation}. However, for the range of problems considered here using \eqref{equ: gap energy minimisation} has delivered both stable and accurate results.

The method to solve this problem is presented in Algorithm \ref{alg: gap func min}. This will update the GIMP position from its uncorrected state, Figure \ref{fig: GIMP contact reset bad}, to its new state, Figure \ref{fig: GIMP contact reset good}. For this algorithm an explicit forward Euler method was used. An implicit method was considered, however during numerical testing it was observed that during the solve it moved the GIMP so it was out of contact, meaning $g_{N,n}^{new}=0~\forall n$ causing the convergence of Algorithm \ref{alg: gap func min} to stagnate. In this regard the forward Euler approach was considered more robust, with the step size of the method constrained by the minimum side length, $h_p$, of the GIMP in contact, divided by $200$.

\begin{algorithm}[H]
\SetAlgoLined
    Calculate $E^{old} =  \sum_{n}^8\left(\epsilon_{N,n}^{old}g_{N,n}^{old}\right)^2$\;
    \texttt{exit} $\leftarrow 0$\;
    \While{\texttt{exit} $= 0$}{
        Calculate the residual $R_i = \frac{\partial E}{\partial x_i}$ \;
        \eIf{$ E /E^{old} > 0.01$}{
            Calculate update step vector $v_i = \frac{R_ih_p}{200|R_i|}$ \;
            Update the MP position $x_i \leftarrow x_i - v_i$ \;
            Recalculate $E$ \;
        }
       {
            \texttt{exit} $\leftarrow 1$\;
        }
    }
\caption{Energy Minimisation Algorithm.}\label{alg: gap func min}
\end{algorithm}
The factor of $h_p/200$ was chosen as compromise, as the algorithm is entirely multiplicative and fast, however reducing the step size beyond a reasonable value could make it expensive. $h_p/200$ is also an upper bound to ensure the algorithm converges and to limit the possibility of the GIMP being pushed out contact.

\section{Numerical simulations}
This section will demonstrate the capabilities of the proposed method via four numerical examples.  The first two examples test the method against analytical solutions to check specific aspects of the contact formulation. The third and fourth examples are comparisons with physical modelling results to demonstrate that the method can obtain robust, physically correct results without parameter tuning.

\subsection{Cube under compression}
\subsubsection*{Example scope}
A 1D contact problem with an analytical solution, but modelled in 3D, is used to validate the normal contact formulation. The validation will investigate the stress solution of one of the GIMPs in contact with the rigid body, as an inconsistent formulation will show poor agreement with the analytical solution \cite{bird2023cone}. Also considered is the stress solution and the global error when the normal contact penalty is varied. In the limit as the penalty$\rightarrow\infty$ the error in the Signorini-Hertz-Moreau conditions should tend to zero and hence the stress solution converge to the true solution.
 
\subsubsection*{Setup}
\begin{figure}[ht!]
    \centering
     \def\svgwidth{0.5\textwidth}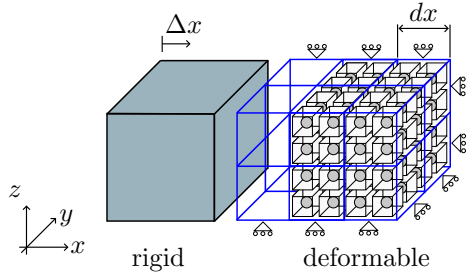
    \caption{Cube under compression: The initial geometry and mesh setup.}
    \label{fig:compression 1D}
\end{figure}

The geometry of the problem is shown in Figure \ref{fig:compression 1D} and consists of two unit sized cubes, one rigid and one deformable. The rigid body is displaced in the $x$-direction by $\Delta x$ whilst the deformable body has roller boundary conditions on all mesh boundaries. As shown in Figure \ref{fig:compression 1D}, the GIMP faces which will be in contact with the rigid body are not consistent with a mesh boundary, hence are free to move subject to the contact.

The deformable body is elastic with Young's Modulus $E = 10^3$Pa, Poisson's ratio $\nu = 0$ and a contact normal penalty of $\epsilon_N = p_fE$Pa, where $p_f$ is the penalty factor that is used to very the contact penalty. It is discretised uniformly with a Cartesian mesh with element edge length $dx=0.1$m. Each element initially contains 2 GIMPs in each direction with edge length $\frac{1}{2}dx$.

Initially, the distance between the rigid and deformable body is $10^{-3}$m, the rigid body is then displaced $\Delta z -= 0.2001$m over $5$ time steps. The error in the $L^2$ norm is presented as
\begin{equation}\label{eq: L2 stress error}
    e = \left({ \sum_{\forall p} |\sigma - \sigma_{p}|^2 V^0_{p} } \right)^{1/2} \quad\text{with}\quad \sigma = E \log\left(\frac{l}{l_0}\right)\left(\frac{l_0}{l}\right)
\end{equation}
as the analytical stress solution, which is constant, where $l$ is the deformed length of the column, measured here using the final position of the rigid body contact surface, and $l_0$ is the original length of the deformable body with length $1$m. $V^0_{p}$ the original volume of the GIMP, and $\sigma_{p}$ is the stress at the GIMP.

\subsubsection*{Results discussion}
\begin{figure}
\begin{subfigure}{0.48\linewidth}
        \centering
        \includegraphics[width=\linewidth]{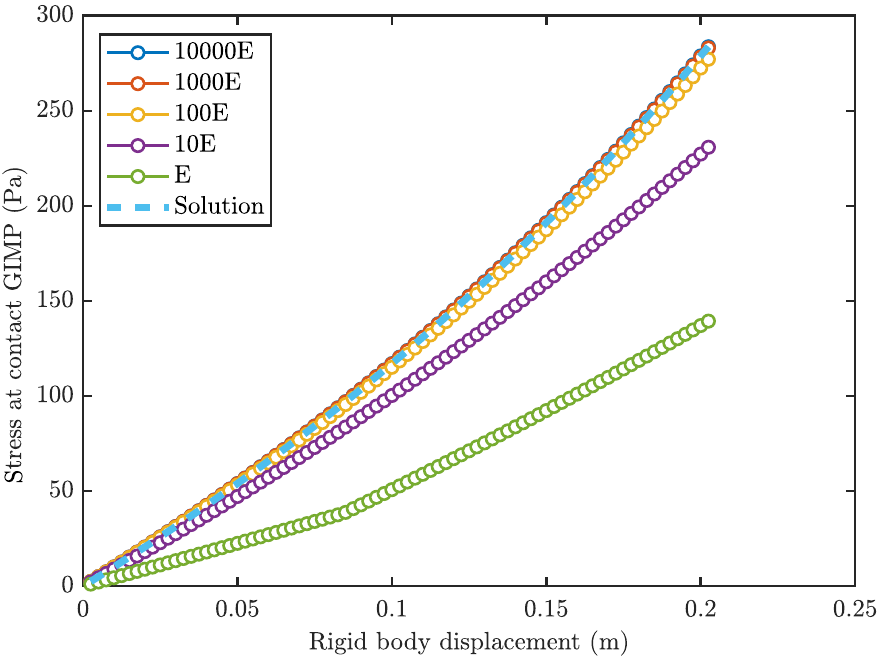}
        \caption{}
        \label{fig:compression_1D_result1}
    \end{subfigure}
    \hfill
    \begin{subfigure}{0.48\linewidth}
        \centering
        \includegraphics[width=\linewidth]{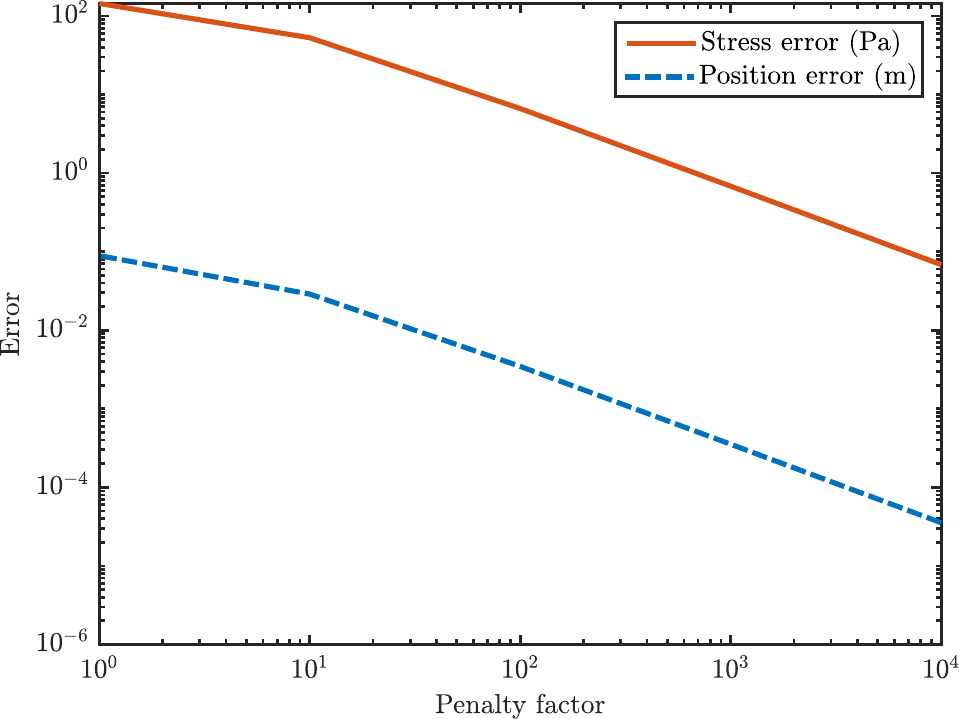}
        \caption{}
        \label{fig:compression_1D_result2}
    \end{subfigure}
    \caption{Cube under compression: (a) is the stress at the contact GIMP for different penalty values, with the corresponding error for a rigid body displacement of $0.2$m.}
    \label{}
\end{figure}

The results presented here demonstrate that the stress at GIMPs in contact with the rigid body is consistent and converges to the analytical solution when the penalty is increased for a range of rigid body displacements, as shown in Figure \ref{fig:compression_1D_result1}. The results indicate that when $p_f$ is near $100E$ there is good agreement between the numerical and analytical results for this problem. Additionally, Figure \ref{fig:compression_1D_result2} shows the global stress solution error and the position error, the interpenetration between the GIMP in contact and the rigid body. Both error measures decrease as the penalty factor increases. These figures therefore show that the stress solution for the GIMP in contact with the rigid body is consistent and that when the penalty is increased the overlap between the two bodies converges to zero.

\subsection{Rolling sphere}
\subsubsection*{Example scope}
The previous problem demonstrated optimal convergence for a pseudo-static problem which considered only normal contact. Here the normal and frictional contact formulation is validated by modelling a sphere rolling down a slope. The numerical results are compared to the problem's analytical solution for a range of friction coefficients, testing the numerical framework's ability to accurately capture the slip and stick states of Coulomb friction.

\subsubsection*{Setup}
\begin{figure}[ht!]
    \centering
    \begin{subfigure}[b]{0.65\textwidth}
        \centering
        \def\svgwidth{\textwidth}
\begingroup%
  \makeatletter%
  \providecommand\color[2][]{%
    \errmessage{(Inkscape) Color is used for the text in Inkscape, but the package 'color.sty' is not loaded}%
    \renewcommand\color[2][]{}%
  }%
  \providecommand\transparent[1]{%
    \errmessage{(Inkscape) Transparency is used (non-zero) for the text in Inkscape, but the package 'transparent.sty' is not loaded}%
    \renewcommand\transparent[1]{}%
  }%
  \providecommand\rotatebox[2]{#2}%
  \newcommand*\fsize{\dimexpr\f@size pt\relax}%
  \newcommand*\lineheight[1]{\fontsize{\fsize}{#1\fsize}\selectfont}%
  \ifx\svgwidth\undefined%
    \setlength{\unitlength}{459.28366353bp}%
    \ifx\svgscale\undefined%
      \relax%
    \else%
      \setlength{\unitlength}{\unitlength * \real{\svgscale}}%
    \fi%
  \else%
    \setlength{\unitlength}{\svgwidth}%
  \fi%
  \global\let\svgwidth\undefined%
  \global\let\svgscale\undefined%
  \makeatother%
  \begin{picture}(1,0.5011326)%
    \lineheight{1}%
    \setlength\tabcolsep{0pt}%
    \put(0,0){\includegraphics[width=\unitlength,page=1]{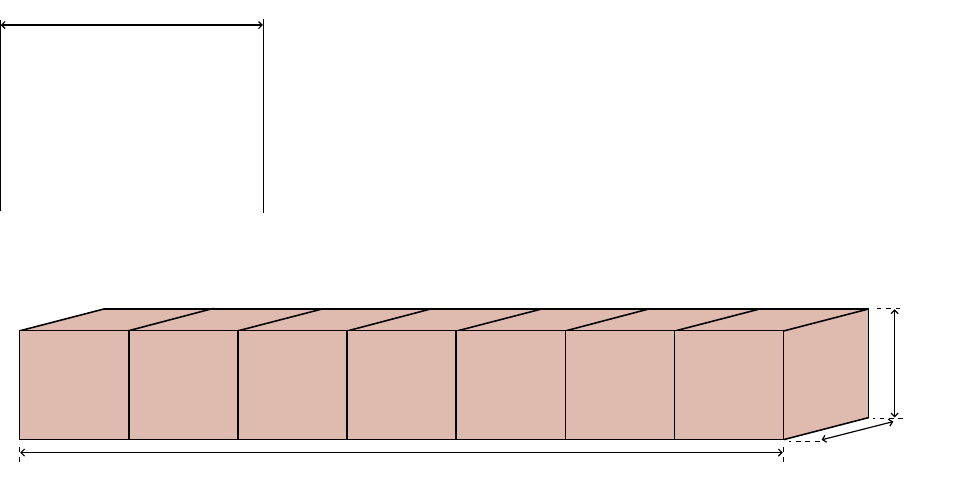}}%
    \put(0.11486656,0.48459023){\makebox(0,0)[lt]{\lineheight{1.25}\smash{\begin{tabular}[t]{l}$d_b$\end{tabular}}}}%
    \put(0.93965649,0.11290442){\makebox(0,0)[lt]{\lineheight{1.25}\smash{\begin{tabular}[t]{l}$L_z$\end{tabular}}}}%
    \put(0.41446017,0.00513494){\makebox(0,0)[lt]{\lineheight{1.25}\smash{\begin{tabular}[t]{l}$L_x$\end{tabular}}}}%
    \put(0.89560082,0.03041108){\makebox(0,0)[lt]{\lineheight{1.25}\smash{\begin{tabular}[t]{l}$L_y$\end{tabular}}}}%
    \put(0,0){\includegraphics[width=\unitlength,page=2]{Rolling_sphere.pdf}}%
    \put(0.15399987,0.22283815){\makebox(0,0)[lt]{\lineheight{1.25}\smash{\begin{tabular}[t]{l}$g_i$\end{tabular}}}}%
    \put(0,0){\includegraphics[width=\unitlength,page=3]{Rolling_sphere.pdf}}%
  \end{picture}%
\endgroup%

        \caption{}
        \label{fig:Sphere_setup}
    \end{subfigure}
    \hfill
    \begin{subfigure}[b]{0.34\textwidth}
        \centering
        \def\svgwidth{\textwidth}
\begingroup%
  \makeatletter%
  \providecommand\color[2][]{%
    \errmessage{(Inkscape) Color is used for the text in Inkscape, but the package 'color.sty' is not loaded}%
    \renewcommand\color[2][]{}%
  }%
  \providecommand\transparent[1]{%
    \errmessage{(Inkscape) Transparency is used (non-zero) for the text in Inkscape, but the package 'transparent.sty' is not loaded}%
    \renewcommand\transparent[1]{}%
  }%
  \providecommand\rotatebox[2]{#2}%
  \newcommand*\fsize{\dimexpr\f@size pt\relax}%
  \newcommand*\lineheight[1]{\fontsize{\fsize}{#1\fsize}\selectfont}%
  \ifx\svgwidth\undefined%
    \setlength{\unitlength}{191.24554491bp}%
    \ifx\svgscale\undefined%
      \relax%
    \else%
      \setlength{\unitlength}{\unitlength * \real{\svgscale}}%
    \fi%
  \else%
    \setlength{\unitlength}{\svgwidth}%
  \fi%
  \global\let\svgwidth\undefined%
  \global\let\svgscale\undefined%
  \makeatother%
  \begin{picture}(1,0.82255737)%
    \lineheight{1}%
    \setlength\tabcolsep{0pt}%
    \put(0,0){\includegraphics[width=\unitlength,page=1]{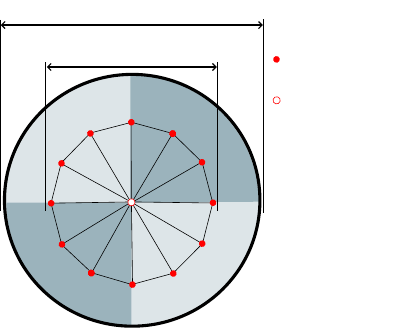}}%
    \put(0.27585631,0.78283025){\makebox(0,0)[lt]{\lineheight{1.25}\smash{\begin{tabular}[t]{l}$d_b$\end{tabular}}}}%
    \put(0.72604605,0.66286383){\makebox(0,0)[lt]{\lineheight{1.25}\smash{\begin{tabular}[t]{l}truss node\end{tabular}}}}%
    \put(0.72772421,0.58475251){\makebox(0,0)[lt]{\lineheight{1.25}\smash{\begin{tabular}[t]{l}zero mass\\\end{tabular}}}}%
    \put(0.27723197,0.68586063){\makebox(0,0)[lt]{\lineheight{1.25}\smash{\begin{tabular}[t]{l}$d_m$\end{tabular}}}}%
    \put(0.72685618,0.53740622){\makebox(0,0)[lt]{\lineheight{1.25}\smash{\begin{tabular}[t]{l}truss node\end{tabular}}}}%
  \end{picture}%
\endgroup%

        \caption{}
        \label{fig:Sphere_truss}
    \end{subfigure}
    \caption{Rolling sphere: Initial geometry and mesh setup in shown in (a) and the corresponding truss frame in (b).}
    \label{fig:Sphere_combined}
\end{figure}

Figure \ref{fig:Sphere_setup} shows the loading, geometry, and mesh for this problem. Although the problem models a sphere rolling down a $45^\circ$ slope, the slope is represented as horizontal to ensure a smooth interface between the sphere and the GIMPs. To simulate the angled slope, gravitational acceleration is applied at an angle corresponding to the original slope, with $g_i = 9.81 \times \left[{1}/{\sqrt{2}}, 0, -{1}/{\sqrt{2}}\right]^\top$m/s$^{2}$.

The kinematics of the rigid body are modelled by a series of elastic truss elements with only displacement degrees of freedom. Since the rigid body is free to rotate and translate, both the translation and rotational inertia need to be represented correctly. Although this is a 3D problem, the sphere is restrained in the $y$-direction and hence the rotational inertia needs only to be correct about the $y$-axis. Therefore the truss elements are arranged as in Figure \ref{fig:Sphere_truss}. At the centre of the sphere is a node with zero mass, which is surrounded by a ring of equally spaced nodes on the y-plane which do have mass. To ensure the rotational inertia of the ring is consistent with the sphere's rotational inertia the diameter of the ring is set to $d_m = d_b\sqrt{(2/5)}$. Here $100$ nodes are used each having a mass of $1/100$kg.

The sphere has a diameter $d_b = 1$m and is constructed from 3120 triangles arranged on a latitude-longitude grid. The smallest and most slender triangles are at the poles whilst the most regular and largest triangles are at the latitude $0^\circ$. The material domain has dimensions, $L_x = 10$m,  $L_y = 1$m,  $L_z = 1$m. The elements in the domain are of unit size and individually contain 4 material points in each direction.

The analytical solution for this problem is defined for two rigid bodies in contact. In order for the deformable body to behave like a rigid body, homogeneous Dirichlet BCs were applied on all its exterior surfaces. The normal and tangential contact penalties for this problem are set as $50E$ and $25E$, with $E=10^5$Pa.

The analytical solution to this problem is presented as the distance $d_x$ the sphere has moved down the slope as a function of time. The solution has two forms depending if the  frictional contact is in a slip or stick state
\begin{equation}
    d_x(t) = \left\{
    \begin{tabular}{ll}
         $(gt^2/2)\left[\sin(\theta_s)-\mu\cos(\theta_s)\right]$ & if slipping, $\tan(\theta_s)>3\mu$  \\
         $gt^2\sin(\theta_s)/3$ & else sticking
    \end{tabular}\right.
\end{equation}
where $g=|g_i|$ is the gravitational acceleration, $\theta_s = 45^\circ$ is the slope angle, $t$s is time and $\mu\in\{0,0.1,0.2,0.4,1.0\}$ are the friction coefficients  considered here.

\subsubsection*{Results discussion}
\begin{figure}
    \centering
    \includegraphics[width=0.5\linewidth]{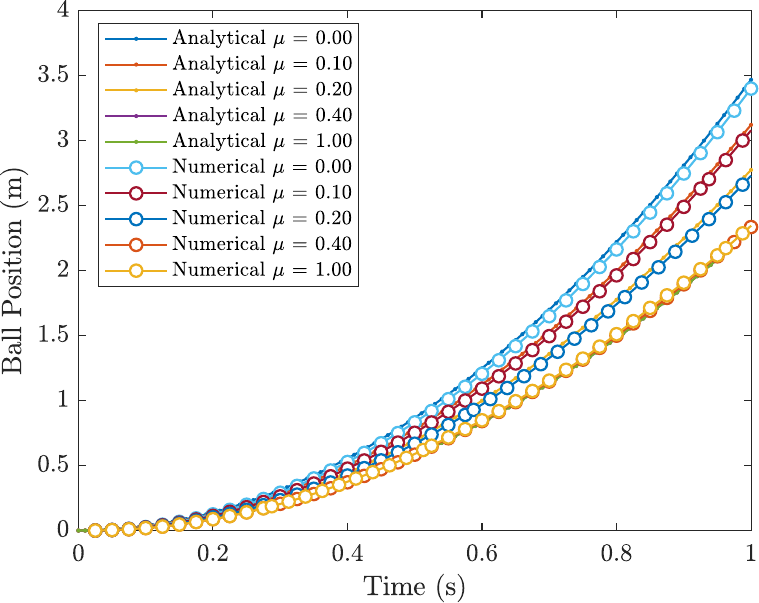}
  \caption{Rolling sphere: The position of the sphere as a function of time for different coefficients of friction.}
    \label{fig:rolling sphere results}
\end{figure}

A comparison of the numerical and analytical results is shown in Figure \ref{fig:rolling sphere results}.
For $\mu\in\{0,0.1,0.2\}$ the sphere is in a slip state, with increasing $\mu$, less slip occurs and hence the slower ball velocity. Over this range Figure \ref{fig:rolling sphere results} shows good agreement, with the numerical results showing a very slight under-prediction of the ball position with time. Good agreement is also seen for the stick states, $\mu\in\{0.4,1.0\}$, demonstrating that once the sphere is sticking the solution is largely invariant to $\mu$.

\subsection{Cone penetration test}
\subsubsection*{Example scope}
In this section a pseudo-static Cone Penetration Test (CPT) into dry sand with a homogeneous relative density is performed. The cone resistance, the vertical force acting on the cone, is then compared to experimental results obtained by \citet{davidson2022physical} for two different relative densities of sand, and for a range of mesh refinements. 

\subsubsection*{Setup}
\begin{figure}[ht!]
    \centering
    \begin{subfigure}[b]{0.48\textwidth}
        \centering
       \def\svgwidth{\textwidth}
\begingroup%
  \makeatletter%
  \providecommand\color[2][]{%
    \errmessage{(Inkscape) Color is used for the text in Inkscape, but the package 'color.sty' is not loaded}%
    \renewcommand\color[2][]{}%
  }%
  \providecommand\transparent[1]{%
    \errmessage{(Inkscape) Transparency is used (non-zero) for the text in Inkscape, but the package 'transparent.sty' is not loaded}%
    \renewcommand\transparent[1]{}%
  }%
  \providecommand\rotatebox[2]{#2}%
  \newcommand*\fsize{\dimexpr\f@size pt\relax}%
  \newcommand*\lineheight[1]{\fontsize{\fsize}{#1\fsize}\selectfont}%
  \ifx\svgwidth\undefined%
    \setlength{\unitlength}{232.04409622bp}%
    \ifx\svgscale\undefined%
      \relax%
    \else%
      \setlength{\unitlength}{\unitlength * \real{\svgscale}}%
    \fi%
  \else%
    \setlength{\unitlength}{\svgwidth}%
  \fi%
  \global\let\svgwidth\undefined%
  \global\let\svgscale\undefined%
  \makeatother%
  \begin{picture}(1,1.24694698)%
    \lineheight{1}%
    \setlength\tabcolsep{0pt}%
    \put(0,0){\includegraphics[width=\unitlength,page=1]{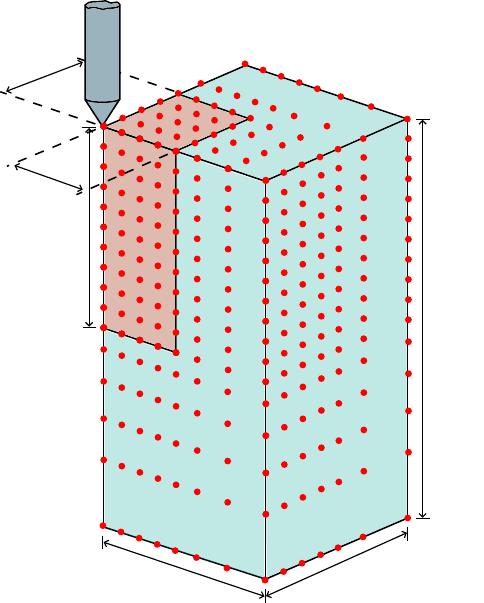}}%
    \put(0.32514366,0.03256806){\makebox(0,0)[lt]{\lineheight{1.25}\smash{\begin{tabular}[t]{l}$L_x$\end{tabular}}}}%
    \put(0.70404636,0.04154104){\makebox(0,0)[lt]{\lineheight{1.25}\smash{\begin{tabular}[t]{l}$L_y$\end{tabular}}}}%
    \put(0.88056246,0.57897493){\makebox(0,0)[lt]{\lineheight{1.25}\smash{\begin{tabular}[t]{l}$L_z$\end{tabular}}}}%
    \put(0.05523343,0.84379916){\makebox(0,0)[lt]{\lineheight{1.25}\smash{\begin{tabular}[t]{l}$U_x$\end{tabular}}}}%
    \put(0.03948152,1.10112002){\makebox(0,0)[lt]{\lineheight{1.25}\smash{\begin{tabular}[t]{l}$U_y$\end{tabular}}}}%
    \put(0.11825199,0.77542279){\makebox(0,0)[lt]{\lineheight{1.25}\smash{\begin{tabular}[t]{l}$U_z$\end{tabular}}}}%
  \end{picture}%
\endgroup%
 
        \caption{}
       \label{fig:CPT setup}
    \end{subfigure}
    \hspace{0.5cm}
    \begin{subfigure}[b]{0.38\textwidth}
        \centering
        \includegraphics[width=\textwidth]{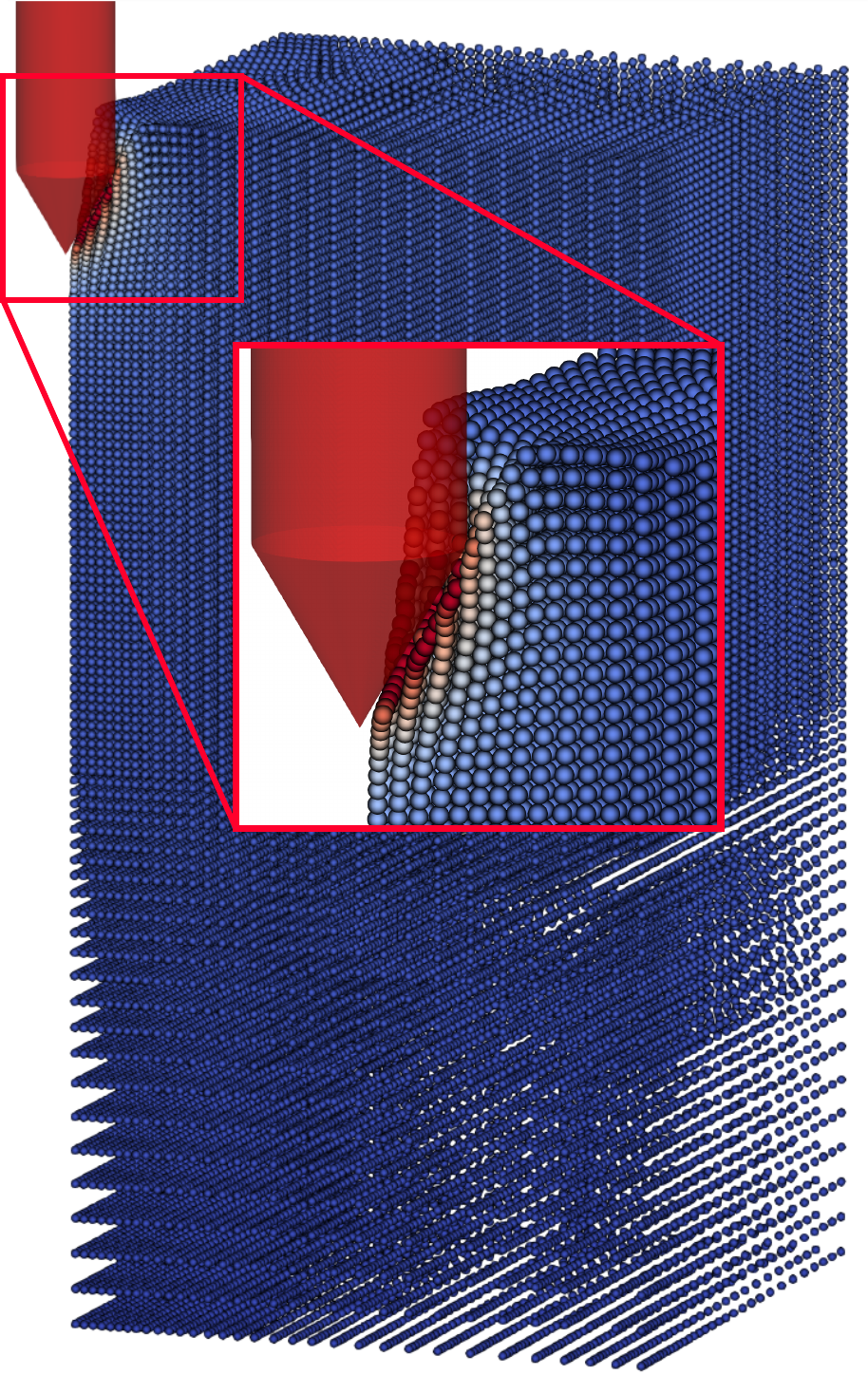}
        \caption{}
        \label{fig:CPT result setup}
    \end{subfigure}
    \caption{Cone penetration test: The initial geometry and mesh setup is shown in (a) with the CPT penetrated $1.2$m shown in (b). The displacement magnitude is shown in (b), blue is $0$m and red is $0.5$m.}
\end{figure}
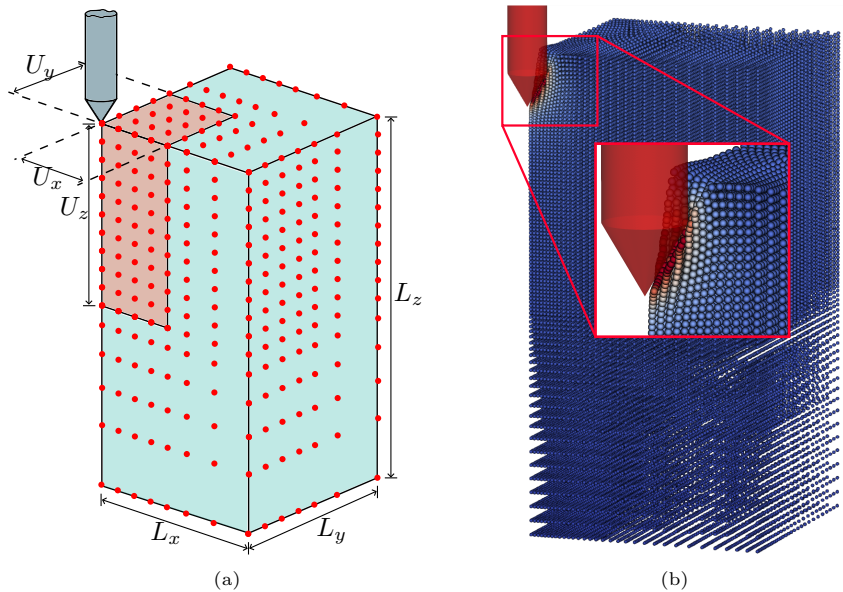

The geometry and node distribution of the problem is shown in Figure \ref{fig:CPT setup}. Provided that the external boundaries are sufficiently far from the cone, the stress solution around the cone is symmetric in the circumferential direction. Problem symmetry could be exploited to reduce computational cost, as in the axisymmetric model of \citet{bird2024implicit}, however in this case one quarter of the problem is analysed to test the 3D implementation. This is done by setting roller boundary conditions at $x=0$m and $y=0$m that represent planes of symmetry. On the remaining external boundaries, $x = L_x$m, $y = L_y$m and $z = 0$m, and roller boundary conditions are also applied. The cone's displacements are fully prescribed; in total the CPT is displaced $-5$m in the $z$-direction over 200 equal time steps, with equal load increment.  However, if a time step fails the step size is reduced by a factor of 2 until convergence is achieved.  The step size is then reset to the original size and the analysis progressed until the specified displacement is achieved.

\begin{table}[ht!]
\centering
\caption{Cone penetration test: deformable domain lengths.}\label{tab: CPT domain sizes}
\label{tab:material}
\begin{tabular}{|r|c|c|c|c|c|c|}
\hline
\textbf{Dimension} & $L_x$ & $L_y$ & $L_z$ & $U_x$ & $U_y$ & $U_z$ \\ \hline
\textbf{Length (m)} &  5  & 5 & 10 & 2.5 & 2.5 & 5 \\
\hline
\end{tabular}
\end{table}

The dimensions of the body are provided in Table \ref{tab: CPT domain sizes}. The red portion of the domain, shown in Figure \ref{fig:CPT setup}, is defined by the dimensions $U_x$, $U_y$ and $U_z$. In this region the discretisation is uniform with length $dx$. In the remainder of the domain, the blue region, the distances between nodes are scaled with a power law from the boundary of the red region to the end of the domain,
\begin{equation}
    dx_{i+1} = (dx_{i})^{1.3} 
\end{equation}
where $dx_{i}$ is the distance between the node $i$ and node $i-1$ in a particular direction. An indication of this nodal distribution is outlined in Figure \ref{fig:CPT setup} with the red dots.

The material used in the real experiment was dry silica sand obtained from Congleton in the UK \cite{davidson2022physical}. There are many sophisticated models for modelling sand however from the justifications and results observed in \citet{bird2024implicit} for modelling CPTs, a linear-elastic, perfectly-plastic material with a Drucker-Prager yield surface is sufficient. All the properties for the Drucker-Prager yield surface can be determined using the equations provided by \citet{brinkgreve2010validation}, and for an overview of the use of the equations see \citet{bird2024implicit}. However the initial state of the material due to gravity also has to be considered and hence the Young's modulus will vary with depth. If a negligible cohesion is assumed (here we use $c = 300$Pa), the variation of the Young's modulus can be determined using the formula provided by \citet{schanz2019hardening}
\begin{equation}\label{eq: E}
    E_{50} = {E}_{50}^{ref}\left(\frac{{{\sigma}_{v} {K}_{0}}}{{{p}^{ref}}} \right)^{m_E}
    \qquad\text{with}\qquad 
    \sigma_v = d_p\rho,
\end{equation}
where $K_0=1-\sin(\phi)$ is the coefficient of earth pressure at rest \cite{jaky1944coefficient}, $\sigma_v$ is the vertical stress and $d_p$ is the distance of the MP from the surface of the sample (i.e. the depth), $m_E$ is an exponent controlling the variation of stiffness \cite{brinkgreve2010validation}. Here two relative densities are considered, $38\%$ and $82\%$ with the values of their material properties provided in Table \ref{tab: material properties}, alongside $44\%$ which is used in Section~\ref{sec:plough}. 

\begin{table}[ht!]
\centering
\caption{Cone penetration test: Material properties.}\label{tab: material properties}
\begin{tabular}{|l|c|c|c|}
\hline
\textbf{Property} & \textbf{38\%} & \textbf{44\%}  &\textbf{82\%} \\
\hline
Reference Young's modulus, $E^{ref}_{50}$ (kPa) & 22,800 & 26,400 & 49,200 \\
Density, $\rho$ (kN/m$^3$)                      & 16.5   & 16.7   & 18.2 \\
Poisson's ratio                                 & 0.3    & 0.3    & 0.3 \\
Friction angle ($^\circ$)                       & 32.8   & 33.5   & 38.3 \\
Dilation angle ($^\circ$)                       & 2.8    & 3.5    & 8.3 \\
Apparent cohesion (kPa)                         & 0.3    & 0.3    & 0.3 \\
{\color{black}Coefficient of earth pressure at rest}, $K_0$ & 0.45 & 0.44 & 0.38 \\
Stiffness exponent, $m_E$                                   & 0.58 & 0.57 & 0.44 \\
\hline
\end{tabular}
\end{table}

The normal and tangential contact penalties are respectively $50E_p$ and $25E_p$, where $p$ is the GIMP that is contact with the material point. This variation in the penalties is required since the material has an inhomogeneous Young's modulus, and the variation enables a more accurate representation of the conduct boundary conditions without making the non-linear solver unstable through an excessively high penalty.

\subsubsection*{Results discussion}
\begin{figure}[ht!]
    \centering
    \includegraphics[width=0.54\textwidth]{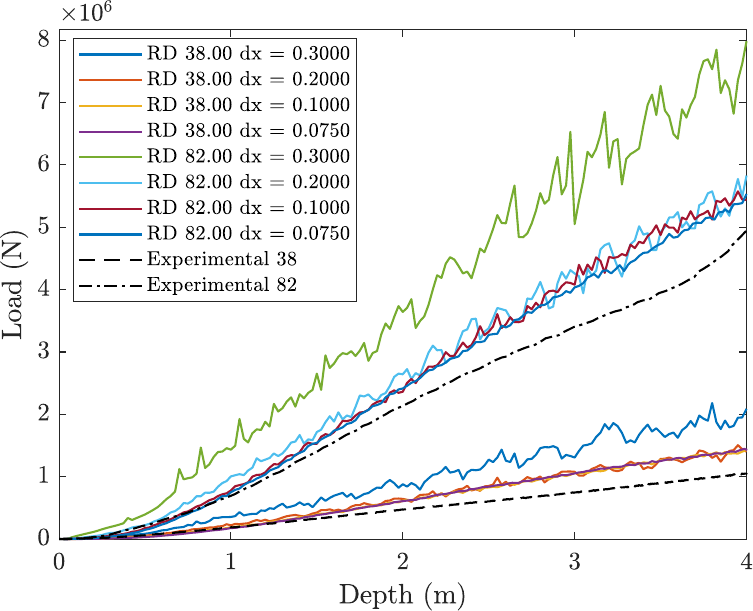}
    \caption{Cone penetration test: Results for CPT tip load with depth.} \label{fig:CPT results}
\end{figure}

The results for the CPT tip load with depth are shown in Figure \ref{fig:CPT results} for a range of refinements, $dx = \{0.075,0.1,0.2,0.3\}$m. 

Figure \ref{fig:CPT results} shows that for the most refined mesh, $dx=0.075$m, good agreement is obtained for both relative densities. It is important to note that these results are pure predictions and there is no tuning of the parameters for the material model, or numerical parameters, performed to obtain a better agreement. It is also observed that with refinement the solutions converge in the same direction towards the experimental result, and in addition the magnitude of the oscillations also decreases. The most likely cause for the oscillations are material points coming in and out of contact with the cone. With refinement each GIMP in contact carries a smaller portion of the load and therefore has a smaller effect on the global response.

 \subsection{Plough}\label{sec:plough}

\subsubsection*{Example scope}
\begin{figure}[ht!]
    \centering
    \def\svgwidth{0.6\textwidth}\input{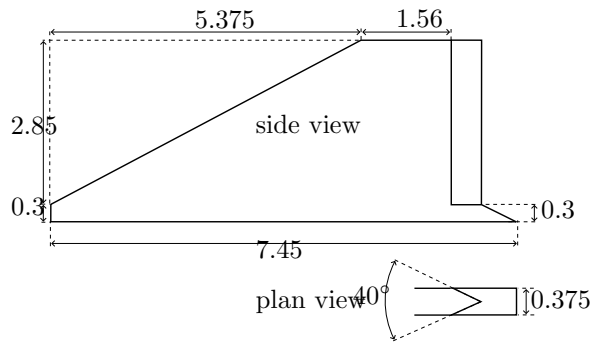}
    \caption{Plough: Dimensioned drawing of the plough, metres.}
    \label{fig:plough dimensions}
\end{figure}

The final problem considered is the modelling of a seabed cable plough, with the geometry shown in Figure \ref{fig:plough dimensions}. The plough is pulled through dry sand at a constant speed and depth. The required horizontal pull (or tow) force predicted by the numerical model is compared to the steady state experimental data obtained by \citet{robinson2021cone} in a $50g$ geotechnical centrifuge. To simplify the numerical analysis, the numerical model considers the full scale problem at $1g$, which the centrifuge is representing. Following the scaling analysis for dry sand by \citet{robinson2019centrifuge} this means that both length and time are scaled equally, hence the plough is pulled the full scale length but at the model speed. The experiment in \cite{robinson2021cone} reaches an equilibrium in both plough depth and pull force. Here the plough will be fixed at the recorded plough depth and the pull force will be compared to the experimental steady state force.

\subsubsection*{Setup}
\begin{figure}[ht!]
    \centering
    \def\svgwidth{0.8\textwidth}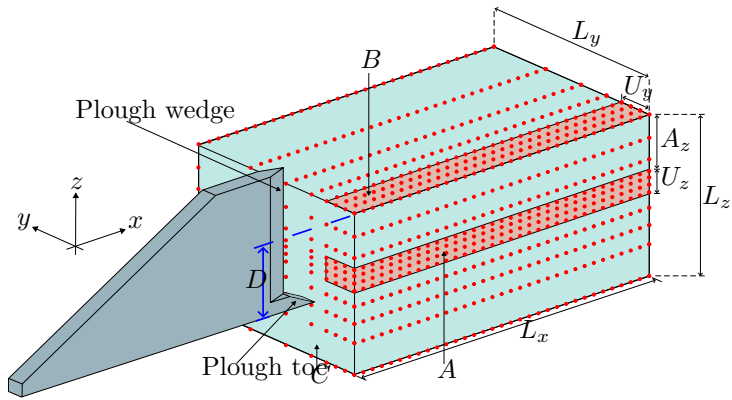
  \caption{Plough: Initial numerical setup of a plough being pulled through dry sand.}
    \label{fig:plough setup}
\end{figure}

The numerical setup for modelling the full plough is shown in Figure \ref{fig:plough setup} with the geometry of the plough provided in Figure \ref{fig:plough dimensions} and corresponding dimension values in Table \ref{tab:plough dimensions}.
\begin{table}[ht!]
\centering
\caption{Plough: deformable domain length.}\label{tab: plough domain sizes}
\label{tab:plough dimensions}
\begin{tabular}{|r|c|c|c|c|c|c|c|}
\hline
\textbf{Dimension} & $L_x$ & $L_y$ & $L_z$ & $U_y$ & $A_z$ & $U_z$& $D$ \\ \hline
\textbf{Length (m)} & 20 & 10 & 7.5 & 1.875 & 1.725 & 0.25& 1.85  \\
\hline
\end{tabular}
\end{table}

The distribution of the initial Cartesian mesh is outlined in Figure \ref{fig:plough setup}, where the red dots show the initial distribution of the element nodes on the boundary of the mesh. To reduce the computational cost only half the problem is modelled, achieved used symmetric boundary conditions in $x$. The purpose of the initial mesh is primarily to setup the initial material point distribution. As the plough travels through the domain, the mesh, but not the material points, is adaptively refined and recovered to reduce the cost associated with the linear solve.
In the red regions the nodes are uniformly distributed with a spacing of $dx$m. In the blue region the node spacing is not uniform and is subject to the power law
\begin{equation}\label{equ plough power law}
    dx_{x+1} = (dx_{i})^{1.3}
\end{equation}
where $dx_i$ is the distance between node $i$ and $i-1$ and $dx_{1} = dx$. If a node does land exactly on the boundary of the domain then the nearest node created by the power law series is snapped to the domain boundary. The node distribution is such that the highest density of nodes, and material points, is around the nose of the plough, marked $A$, where the most complicated 3D displacements occur. Elsewhere where the displacement field is more 2D, such as at the main wedge of the plough, there is a greater refinement in the $x$ and $y$ directions compared to the $z$ direction. This is obtained through the refined region marked $B$ in Figure \ref{fig:plough setup}. The material that is not in contact with the plough has only a coarse refinement. Lastly, the sizes of the refined region sizes will depend slightly on the values of $dx$ so that all elements in this region are equally sized. This can be expressed by the following equation
\begin{equation}
    \bar{U} =  \left\lceil \frac{U}{dx} \right\rceil dx
\end{equation}
where $U$ is the desired size of the refined region and $\bar{U}$ is the actual size. The initial number of GIMPs per elements is $8$, with a $2\times 2 \times 2 $ equally spaces configuration.

\begin{figure}[ht!]
    \centering
    \begin{subfigure}[b]{0.31\textwidth}
        \centering
        \def\svgwidth{\textwidth}
\begingroup%
  \makeatletter%
  \providecommand\color[2][]{%
    \errmessage{(Inkscape) Color is used for the text in Inkscape, but the package 'color.sty' is not loaded}%
    \renewcommand\color[2][]{}%
  }%
  \providecommand\transparent[1]{%
    \errmessage{(Inkscape) Transparency is used (non-zero) for the text in Inkscape, but the package 'transparent.sty' is not loaded}%
    \renewcommand\transparent[1]{}%
  }%
  \providecommand\rotatebox[2]{#2}%
  \newcommand*\fsize{\dimexpr\f@size pt\relax}%
  \newcommand*\lineheight[1]{\fontsize{\fsize}{#1\fsize}\selectfont}%
  \ifx\svgwidth\undefined%
    \setlength{\unitlength}{187.97983725bp}%
    \ifx\svgscale\undefined%
      \relax%
    \else%
      \setlength{\unitlength}{\unitlength * \real{\svgscale}}%
    \fi%
  \else%
    \setlength{\unitlength}{\svgwidth}%
  \fi%
  \global\let\svgwidth\undefined%
  \global\let\svgscale\undefined%
  \makeatother%
  \begin{picture}(1,1.1075312)%
    \lineheight{1}%
    \setlength\tabcolsep{0pt}%
    \put(0,0){\includegraphics[width=\unitlength,page=1]{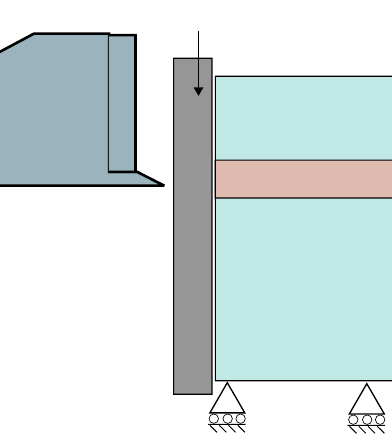}}%
    \put(0.32347364,1.06711389){\makebox(0,0)[lt]{\lineheight{1.25}\smash{\begin{tabular}[t]{l}Frictionless rigid body\end{tabular}}}}%
    \put(0,0){\includegraphics[width=\unitlength,page=2]{plough_side_view.pdf}}%
    \put(0.30560002,0.19945096){\makebox(0,0)[lt]{\lineheight{1.25}\smash{\begin{tabular}[t]{l}$x$\end{tabular}}}}%
    \put(0.05171112,0.4629622){\makebox(0,0)[lt]{\lineheight{1.25}\smash{\begin{tabular}[t]{l}$z$\end{tabular}}}}%
    \put(0,0){\includegraphics[width=\unitlength,page=3]{plough_side_view.pdf}}%
  \end{picture}%
\endgroup%

        \caption{}
        \label{fig:plough side view}
    \end{subfigure}
    \hfill
    \begin{subfigure}[b]{0.57\textwidth}
        \centering
        \def\svgwidth{\textwidth}
\begingroup%
  \makeatletter%
  \providecommand\color[2][]{%
    \errmessage{(Inkscape) Color is used for the text in Inkscape, but the package 'color.sty' is not loaded}%
    \renewcommand\color[2][]{}%
  }%
  \providecommand\transparent[1]{%
    \errmessage{(Inkscape) Transparency is used (non-zero) for the text in Inkscape, but the package 'transparent.sty' is not loaded}%
    \renewcommand\transparent[1]{}%
  }%
  \providecommand\rotatebox[2]{#2}%
  \newcommand*\fsize{\dimexpr\f@size pt\relax}%
  \newcommand*\lineheight[1]{\fontsize{\fsize}{#1\fsize}\selectfont}%
  \ifx\svgwidth\undefined%
    \setlength{\unitlength}{351.10174464bp}%
    \ifx\svgscale\undefined%
      \relax%
    \else%
      \setlength{\unitlength}{\unitlength * \real{\svgscale}}%
    \fi%
  \else%
    \setlength{\unitlength}{\svgwidth}%
  \fi%
  \global\let\svgwidth\undefined%
  \global\let\svgscale\undefined%
  \makeatother%
  \begin{picture}(1,0.5961257)%
    \lineheight{1}%
    \setlength\tabcolsep{0pt}%
    \put(0,0){\includegraphics[width=\unitlength,page=1]{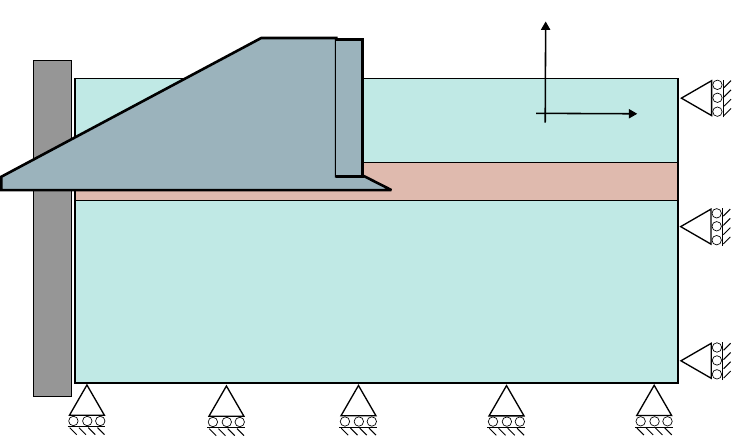}}%
    \put(0.87028187,0.43340239){\makebox(0,0)[lt]{\lineheight{1.25}\smash{\begin{tabular}[t]{l}$x$\end{tabular}}}}%
    \put(0.73434978,0.57448628){\makebox(0,0)[lt]{\lineheight{1.25}\smash{\begin{tabular}[t]{l}$z$\end{tabular}}}}%
    \put(0,0){\includegraphics[width=\unitlength,page=2]{plough_addapt.pdf}}%
    \put(0.46446984,0.26525371){\makebox(0,0)[lt]{\lineheight{1.25}\smash{\begin{tabular}[t]{l}$U_x$\end{tabular}}}}%
    \put(0.40753067,0.14869149){\makebox(0,0)[lt]{\lineheight{1.25}\smash{\begin{tabular}[t]{l}Grid nodes in $x$\end{tabular}}}}%
    \put(0.56761264,0.26535078){\makebox(0,0)[lt]{\lineheight{1.25}\smash{\begin{tabular}[t]{l}$U_x$\end{tabular}}}}%
    \put(0,0){\includegraphics[width=\unitlength,page=3]{plough_addapt.pdf}}%
  \end{picture}%
\endgroup%

        \caption{}
        \label{fig:plough_adaptive_refinement}
    \end{subfigure}
    \caption{Plough: (a) Side view of Signorini boundary conditions with adaptive refinement in $x$ shown in (b).}
\end{figure}

On all faces a roller boundary condition is applied, except on the top face where a homogeneous Neumann boundary exists, and on the face marked $C$, see Figure \ref{fig:plough setup}. On face $C$ a Signorini boundary condition exists so that material can move away from the domain boundary but cannot move across it. Pragmatically, this is enforced with a second rigid body with no friction as shown in Figure \ref{fig:plough side view}.  The plough was incremented at a step size of $0.025$m and there was no interaction between the plough and the secondary rigid body. 

In initial testing, it was observed that the majority of the deformation was occurring at the front the plough, therefore to reduce the computational cost of the simulation mesh adaptation was developed in $x$ as the plough moved through the domain, as shown in Figure \ref{fig:plough_adaptive_refinement}. About the tip of the plough, in both the positive and negative directions, a region with $U_x = 0.5$m was uniformly refined. Outside this region the power law described with Equation \eqref{equ plough power law} was used. To prevent the rigid body penetrating GIMPs all edges with an angle convex edges with an angle less than 90\% were filleted with a radius equal to the minimum side half length of the GIMP, there were approximately 10 segments per $90^\circ$ of fillet.

\subsubsection*{Results discussion}

Here the numerical results for the plough tow force are compared to the equivalent centrifuge experiment, \cite{robinson2019centrifuge}, and additionally the performance of the numerical algorithm is discussed. In the previous sections the geometry of the rigid bodies were relatively simple, and in addition  the point of the CPT coincided with the edge of the domain so there was no risk of the rigid body going inside a material point. For this problem however, the plough contains many surfaces which are connected via rounded concave and convex edges. 

\begin{figure}[ht!]
    \centering
    \includegraphics[width=0.8\linewidth]{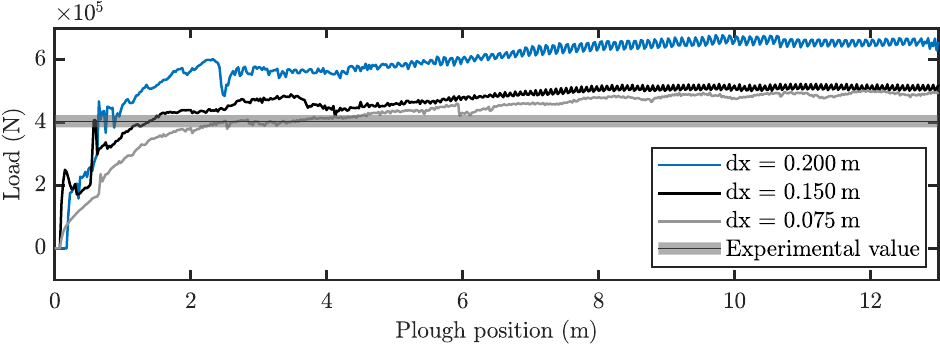}
    \caption{Plough: Comparison of numerical results for pull load with plough position against experimental data for three mesh refinements.}
    \label{fig:plough results}
\end{figure}

\begin{figure}[ht!]
    \centering
    \begin{subfigure}[b]{0.31\textwidth}
        \centering
            \includegraphics[width=\textwidth]{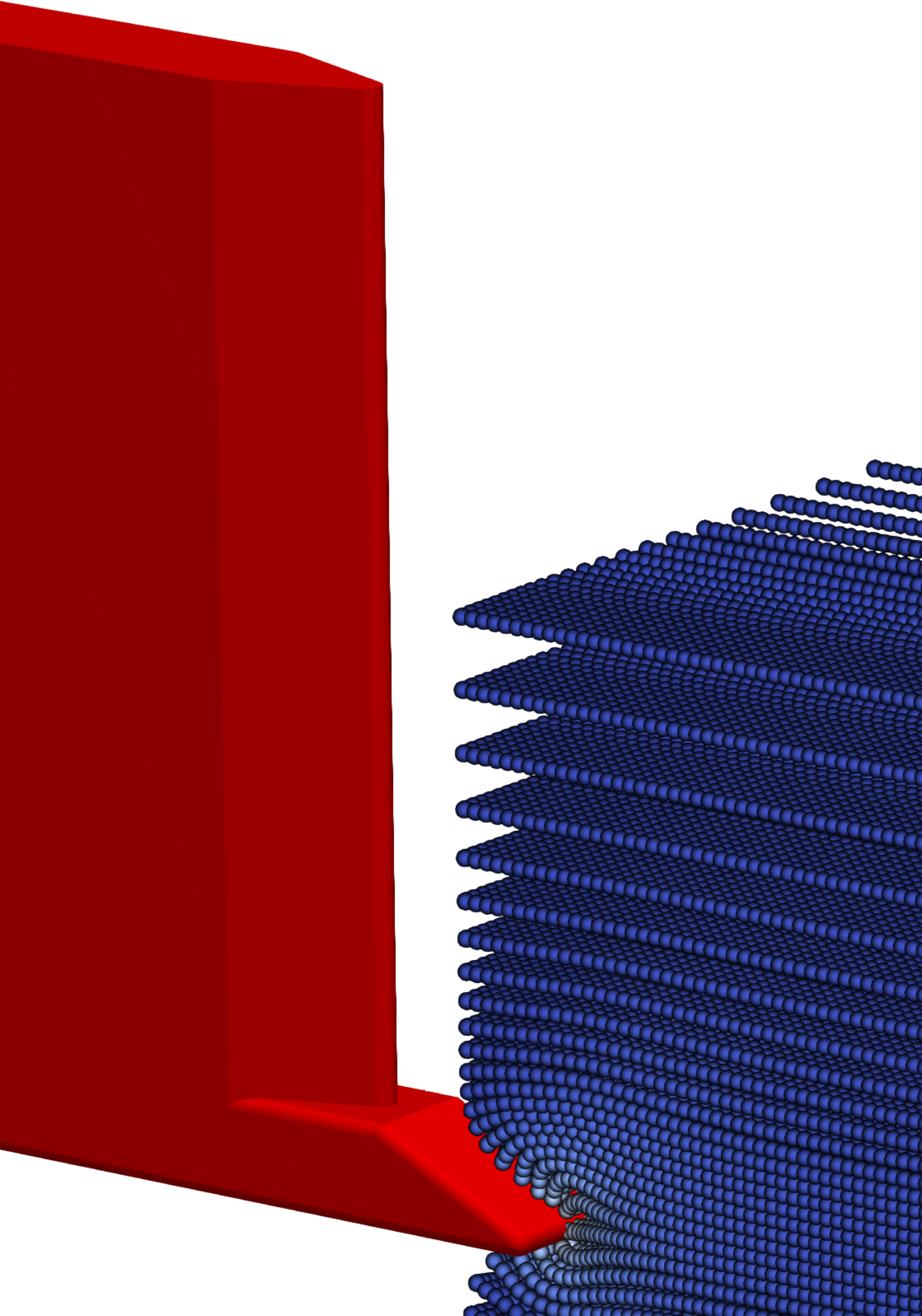}
        \caption{}
        \label{fig:plough before}
    \end{subfigure}
    \hspace{1cm}
    \begin{subfigure}[b]{0.3\textwidth}
        \centering
        \includegraphics[width=\textwidth]{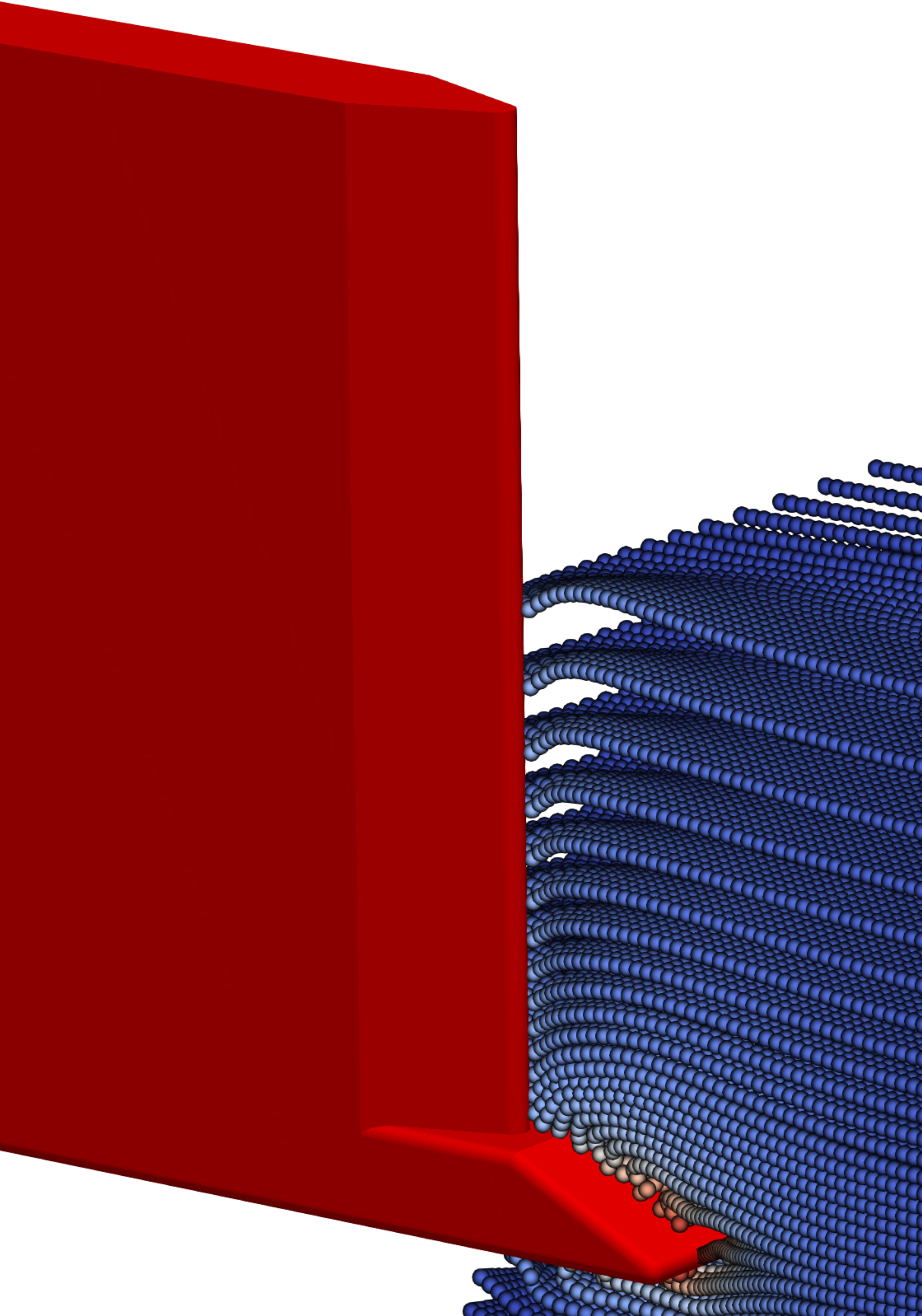}
        \caption{}
        \label{fig:plough after}
    \end{subfigure}
    \caption{Plough: deformed GIMP positions coloured according to the $x$ displacement (red is $1$m blue is $0$m) for a plough just before (a) and just after (b) the front wedge is embedded.}
\end{figure}

Three mesh sizes were considered, $dx = \{0.075,0.15,0.2\}$m and the results of the plough tow force with plough positions are shown in Figure \ref{fig:plough results} alongside the experimental data. The primary observation from the numerical results is that with increasing refinement the tow force gets closer to the experiment, with generally good agreement obtained between mesh sizes $dx\in\{0.075,0.15\}$m and the experimental results. The trend of the numerical results being consistently slightly higher than those observed experimental is consistent with the CPT problem results, and again the simulation is purely predictive and no adjustments to the material or numerical parameters were performed. It is also likely that reducing relative density could close the distance between the results. Before $7.45$m the plough is embedding itself into the domain from the side. In this initial region, $0\rightarrow7.45$m, the results are the most chaotic with sudden changes in the force caused by the build up of material in front of the plough from the initial embedment. Additionally there are sharp peaks in the load caused by the different parts of the plough coming into contact with the domain. For example the greyed region in Figure \ref{fig:plough results} is the plough wedge first coming into contact, going from Figure \ref{fig:plough before} to \ref{fig:plough after}.

\begin{figure}[ht!]
    \centering
    \includegraphics[width=1.0\linewidth]{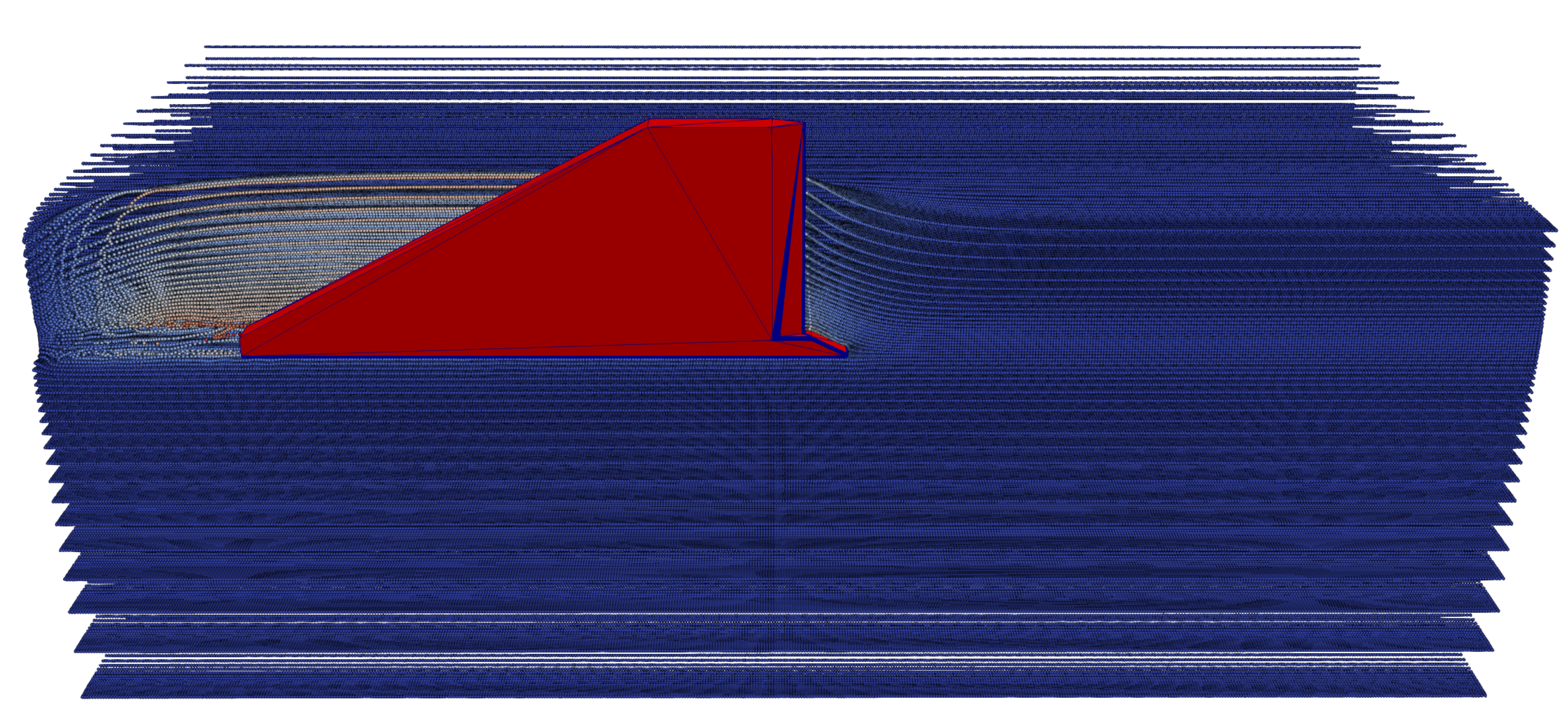}
    \caption{Plough: deformed GIMP positions coloured according to the $x$ displacement (red is $3$m blue is $0$m) for a plough embedded $10$m into the material for the refinement $dx = 0.075$m.}
    \label{fig:plough full 10 m}
\end{figure}

\begin{figure}[ht!]
    \centering
    \includegraphics[width=0.8\linewidth]{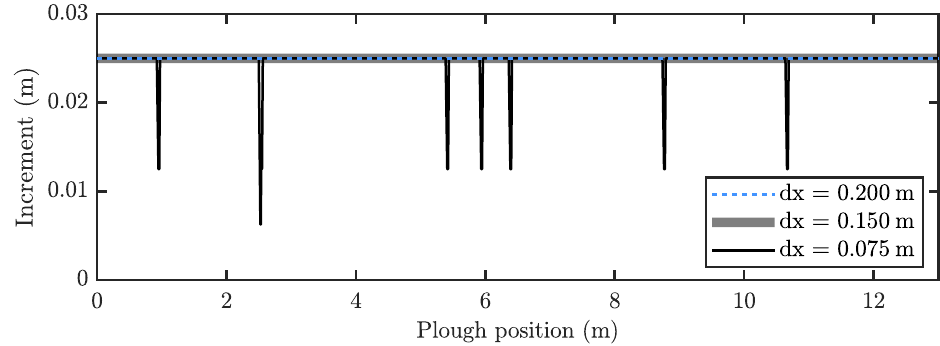}
    \caption{Plough: Comparison of successful load increment sizes for the three mesh sizes.}
    \label{fig:plough results inc}
\end{figure}

The plough is the most complex geometry considered in this paper, and it is therefore appropriate to comment on the stability of this problem. For all meshes the plough was incremented $0.025$m each time step, if the time step failed to converge the step size was reduced by a factor of $2$ and restarted. The displacement increment sizes were recorded for each mesh size and are plotted in Figure \ref{fig:plough results inc}. The figure shows that for all meshes the algorithm remained stable for all mesh sizes, with only $dx=0.075$m requiring a reduction in the displacement step size. This shows that not only good agreement with experimental data was achieved but the algorithms presented here are robust and reliable.

\section{Conclusion}
This paper provides a 3D formulation for contact between a deformable body, modelled using the material point method, and a rigid body. The method uses the definition of the GIMP domain to impose point-to-surface contact between the vertices of the GIMP and the surface of the rigid body. Importantly, no boundary reconstruction is necessary with this approach, and the contact constraints are applied consistently with the numerical representation of the extent of the physical material. It was demonstrated that only when the contact conditions are applied consistently the correct stress field on the contact boundary is obtained (confirming the 2D result of \citet{bird2024implicit}). As the vertices of the GIMP domain can individually be in contact with the rigid body, the possibility of material being undetected by the CPP at rigid body rentrant corner is significantly reduced. This reduces the possibility of material travelling into the rigid body, this is both inaccurate and can cause numerical instability.
%
%
The frictional interaction between the deformable and rigid bodies was validated using the analytical solution of a sphere rolling down a slope. This validated both the normal and the frictional contact, modelled using Coulomb friction, for a rigid body that was free to move. To validate the formulation for examples practical problems, two in the area of soil-structure interaction were chosen with comparisons made against  experimental results obtained from a geotechnical centrifuge, there being no analytical solutions. This included the load on the end of the cone for a CPT and the tow force required to pull a plough at a fixed depth through sand. Both agreed well with experimental data, and both saw convergence towards the experimental result with refinement.  The numerical results were based on a purely predictive methodology, that is, the material constitutive parameters were not tuned to match the experimental results implying that the method can be reliably used to interpolate within and extrapolate beyond the physical dataset. Additionally it was demonstrated for the plough that the solution methodology was stable for difficult contact problems with concave and convex geometries.

\section*{Acknowledgements}
This work was supported by the Engineering and Physical Sciences Research Council [grant numbers EP/W000970/1, EP/W000997/1 and EP/W000954/1]. The second author was supported by funding from the Faculty of Science, Durham University.  All data created during this research are openly available at \href{https://collections.durham.ac.uk/}{collections.durham.ac.uk} (specific DOI to be confirmed if/when the paper is accepted).

\bibliographystyle{plainnat}
\bibliography{main.bib}

\begin{thebibliography}{46}
\providecommand{\natexlab}[1]{#1}
\providecommand{\url}[1]{\texttt{#1}}
\expandafter\ifx\csname urlstyle\endcsname\relax
  \providecommand{\doi}[1]{doi: #1}\else
  \providecommand{\doi}{doi: \begingroup \urlstyle{rm}\Url}\fi

\bibitem[Augarde et~al.(2021)Augarde, Lee, and Loukidis]{Augarde2021}
Charles~E. Augarde, Seung~Jae Lee, and Dimitrios Loukidis.
\newblock Numerical modelling of large deformation problems in geotechnical engineering: A state-of-the-art review.
\newblock \emph{Soils and Foundations}, 61\penalty0 (6):\penalty0 1718--1735, 2021.

\bibitem[Bardenhagen and Kober(2004)]{bardenhagen2004generalized}
S.G. Bardenhagen and E.M. Kober.
\newblock The generalized interpolation material point method.
\newblock \emph{Computer Modeling in Engineering and Sciences}, 5\penalty0 (6):\penalty0 477--496, 2004.

\bibitem[Bing et~al.(2019)Bing, Cortis, Charlton, Coombs, and Augarde]{Bing2019}
Y.~Bing, M.~Cortis, T.J. Charlton, W.M. Coombs, and C.E. Augarde.
\newblock B-spline based boundary conditions in the material point method.
\newblock \emph{Computers \& Structures}, 212:\penalty0 257 -- 274, 2019.

\bibitem[Bird et~al.(2023)Bird, Coombs, Augarde, Brown, Sharif, Carter, Johnson, and Macdonald]{bird2023cone}
Robert Bird, William~M Coombs, Charles Augarde, Michael Brown, Yaseen Sharif, Gareth Carter, Kirstin Johnson, and Catriona Macdonald.
\newblock Cone penetration tests ({CPTs}) in layered soils: a material point approach.
\newblock In \emph{Numerical Methods in Geotechnical Engineering 2023}, 2023.

\bibitem[Bird et~al.(2024)Bird, Pretti, Coombs, Augarde, Sharif, Brown, Carter, Macdonald, and Johnson]{bird2024implicit}
Robert~E Bird, Giuliano Pretti, William~M Coombs, Charles~E Augarde, Yaseen~U Sharif, Michael~J Brown, Gareth Carter, Catriona Macdonald, and Kirstin Johnson.
\newblock An implicit material point-to-rigid body contact approach for large deformation soil--structure interaction.
\newblock \emph{Computers and Geotechnics}, 174:\penalty0 106646, 2024.

\bibitem[Brinkgreve et~al.(2010)Brinkgreve, Engin, and Engin]{brinkgreve2010validation}
R.B.J Brinkgreve, E.~Engin, and H.K. Engin.
\newblock Validation of empirical formulas to derive model parameters for sands.
\newblock \emph{Numerical methods in geotechnical engineering}, 137:\penalty0 142, 2010.

\bibitem[Burman(2010)]{Burman2010}
E~Burman.
\newblock Ghost penalty.
\newblock \emph{Comptes Rendus. Math{\'e}matique}, 348\penalty0 (21):\penalty0 1217--1220, 2010.

\bibitem[Chandra et~al.(2021)Chandra, Singer, Teschemacher, W\"{u}cheer, and Larese]{Chandra2021}
B.~Chandra, V.~Singer, T.~Teschemacher, R.~W\"{u}cheer, and A.~Larese.
\newblock Nonconforming dirichlet boundary conditions in implicit material point method by means of penalty augmentation.
\newblock \emph{Acta Geotechnica}, 16:\penalty0 2315--2335, 2021.

\bibitem[Charlton et~al.(2017)Charlton, Coombs, and Augarde]{charlton2017igimp}
T.J. Charlton, W.M. Coombs, and C.E. Augarde.
\newblock {iGIMP}: An implicit generalised interpolation material point method for large deformations.
\newblock \emph{Computers \& Structures}, 190:\penalty0 108--125, 2017.

\bibitem[Coombs(2023)]{coombs2022ghost}
W.~M. Coombs.
\newblock Ghost stabilisation of the material point method for stable quasi-static and dynamic analysis of large deformation problems.
\newblock \emph{International Journal for Numerical Methods in Engineering}, 124\penalty0 (21):\penalty0 4841--4875, 2023.

\bibitem[Coombs and Augarde(2020)]{coombs2020aample}
W.M. Coombs and C.E. Augarde.
\newblock {AMPLE}: a material point learning environment.
\newblock \emph{Advances in Engineering Software}, 139:\penalty0 102748, 2020.

\bibitem[Coombs et~al.(2020)Coombs, Augarde, Brennan, Brown, Charlton, Knappett, {Ghaffari Motlagh}, and Wang]{Coombs2020continuum}
W.M. Coombs, C.E. Augarde, A.J. Brennan, M.J. Brown, T.J. Charlton, J.A. Knappett, Y.~{Ghaffari Motlagh}, and L.~Wang.
\newblock On {L}agrangian mechanics and the implicit material point method for large deformation elasto-plasticity.
\newblock \emph{Computer Methods in Applied Mechanics and Engineering}, 358:\penalty0 112622, 2020.

\bibitem[Cortis et~al.(2018)Cortis, Coombs, Augarde, Brown, Brennan, and Robinson]{cortis2018imposition}
M.~Cortis, W.M. Coombs, C.E. Augarde, M.~Brown, A.~Brennan, and S.~Robinson.
\newblock Imposition of essential boundary conditions in the material point method.
\newblock \emph{International Journal for Numerical Methods in Engineering}, 113\penalty0 (1):\penalty0 130--152, 2018.

\bibitem[Crisfield(1997)]{Crisfield1997}
M.A. Crisfield.
\newblock \emph{Non-linear Finite Element Analysis of Solids and Structures. Volume 2: Advanced Topics}.
\newblock Wiley, 1997.

\bibitem[Cundall and Strack(1979)]{Cundall1979}
P.~A. Cundall and O.~D.~L. Strack.
\newblock A discrete numerical model for granular assemblies.
\newblock \emph{G\'{e}otechnique}, 29\penalty0 (1):\penalty0 47--65, 1979.

\bibitem[Curnier et~al.(1995)Curnier, He, and Klarbring]{curnier1995continuum}
A.~Curnier, Q.~He, and A.~Klarbring.
\newblock Continuum mechanics modelling of large deformation contact with friction.
\newblock \emph{Contact Mechanics}, pages 145--158, 1995.

\bibitem[Davidson et~al.(2022)Davidson, Brown, Cerfontaine, Knappett, Brennan, Al-Baghdadi, Augarde, Coombs, Wang, Blake, Richards, and Ball]{davidson2022physical}
C.~Davidson, M.~Brown, B.~Cerfontaine, J.~Knappett, A.~Brennan, T.~Al-Baghdadi, C.E. Augarde, W.M. Coombs, L.~Wang, A.~Blake, D.~Richards, and J.D. Ball.
\newblock Physical modelling to demonstrate the feasibility of screw piles for offshore jacket supported wind energy structures.
\newblock \emph{G\'{e}otechnique}, 72\penalty0 (2):\penalty0 108--126, 2022.

\bibitem[de~Souza~Neto et~al.(2008)de~Souza~Neto, Peric, and Owen]{SouzaNeto}
E.~A. de~Souza~Neto, D.~Peric, and D.~R.~J. Owen.
\newblock \emph{Computational Methods For Plasticity: Theory and Applications}.
\newblock John Wiley \& Sons, Ltd, 2008.

\bibitem[de~Vaucorbeil et~al.(2020)de~Vaucorbeil, Nguyen, Sinaie, and Wu]{Vaucorbeil2020}
Alban de~Vaucorbeil, Vinh~Phu Nguyen, Sina Sinaie, and Jian~Ying Wu.
\newblock Material point method after 25 years: Theory, implementation, and applications.
\newblock \emph{Advances in applied mechanics}, 53:\penalty0 185--398, 2020.

\bibitem[Dvorkin and Goldschmit(2006)]{dvorkin2006nonlinear}
Eduardo~N Dvorkin and Marcela~B Goldschmit.
\newblock \emph{Nonlinear continua}.
\newblock Springer Science \& Business Media, 2006.

\bibitem[{Gonzalez Acosta} et~al.(2021){Gonzalez Acosta}, Vardon, and Hicks]{acosta2021development}
J.L. {Gonzalez Acosta}, P.J. Vardon, and M.A. Hicks.
\newblock Development of an implicit contact technique for the material point method.
\newblock \emph{Computers and Geotechnics}, 130:\penalty0 103859, 2021.

\bibitem[Hughes et~al.(1981)Hughes, Liu, and Zimmermann]{Hughes1981}
Thomas~J.R. Hughes, Wing~Kam Liu, and Thomas~K. Zimmermann.
\newblock Lagrangian-eulerian finite element formulation for incompressible viscous flows.
\newblock \emph{Computer Methods in Applied Mechanics and Engineering}, 29\penalty0 (3):\penalty0 329--349, 1981.

\bibitem[Idelsohn et~al.(2004)Idelsohn, O\~{n}ate, and Pin]{Idelsohn2004}
S.R. Idelsohn, E.~O\~{n}ate, and F.~Del Pin.
\newblock The particle finite element method: a powerful tool to solve incompressible flows with free-surfaces and breaking waves.
\newblock \emph{International Journal for Numerical Methods in Engineering}, 61\penalty0 (7):\penalty0 964--989, 2004.

\bibitem[Jaky(1944)]{jaky1944coefficient}
J.~Jaky.
\newblock The coefficient of earth pressure at rest.
\newblock \emph{Journal of the Society of Hungarian Architects and engineers}, 1944.

\bibitem[Lei et~al.(2022)Lei, Wu, Wu, Nie, Cheng, and Zhang]{Lei2022}
Z.~Lei, B.~Wu, S.~Wu, Y.~Nie, S.~Cheng, and C.~Zhang.
\newblock A material point-finite element {(MPM-FEM)} model for simulating three-dimensional soil-structure interactions with the hybrid contact method.
\newblock \emph{Computers and Geotechnics}, 152:\penalty0 105009, 2022.

\bibitem[Ma et~al.(2010)Ma, Giguere, Jayaraman, and Zhang]{Ma2010}
Xia Ma, Paul~T. Giguere, Balaji Jayaraman, and Duan~Z. Zhang.
\newblock Distribution coefficient algorithm for small mass nodes in material point method.
\newblock \emph{Journal of Computational Physics}, 229\penalty0 (20):\penalty0 7819--7833, 2010.

\bibitem[Martinelli and Vahid(2021)]{Martinelli2021}
M.~Martinelli and G.~Vahid.
\newblock Investigation of the material point method in the simulation of cone penetration tests in dry sand.
\newblock \emph{Computers and Geotechnics}, 130:\penalty0 103923, 2021.

\bibitem[Nakamura et~al.(2021)Nakamura, Matsumura, and Mizutani]{nakamura2021particle}
K.~Nakamura, S.~Matsumura, and T.~Mizutani.
\newblock Particle-to-surface frictional contact algorithm for material point method using weighted least squares.
\newblock \emph{Computers and Geotechnics}, 134:\penalty0 104069, 2021.

\bibitem[Noh(1963)]{Noh1963}
W~F Noh.
\newblock {CEL:} a time-dependent, two-space-dimensional, coupled {Eulerian-Lagrange} code.
\newblock 8 1963.
\newblock \doi{10.2172/4621975}.

\bibitem[Pietrzak and Curnier(1999)]{pietrzak1999large}
G.~Pietrzak and A.~Curnier.
\newblock Large deformation frictional contact mechanics: Continuum formulation and augmented {L}agrangian treatment.
\newblock \emph{Computer Methods in Applied Mechanics and Engineering}, 177\penalty0 (3-4):\penalty0 351--381, 1999.

\bibitem[Poulios and Renard(2015)]{poulios2015unconstrained}
Konstantinos Poulios and Yves Renard.
\newblock An unconstrained integral approximation of large sliding frictional contact between deformable solids.
\newblock \emph{Computers \& Structures}, 153:\penalty0 75--90, 2015.

\bibitem[Pretti(2024)]{pretti2024continuum}
G.~Pretti.
\newblock \emph{Continuum mechanics and implicit material point method to underpin the modelling of drag anchors for cable risk assessment}.
\newblock PhD thesis, Durham University, 2024.

\bibitem[Robinson et~al.(2021)Robinson, Brown, Matsui, Brennan, Augarde, Coombs, and Cortis]{robinson2021cone}
S.~Robinson, M.J. Brown, H.~Matsui, A.~Brennan, C.E. Augarde, W.M. Coombs, and M.~Cortis.
\newblock A cone penetration test {(CPT)} approach to cable plough performance prediction based upon centrifuge model testing.
\newblock \emph{Canadian Geotechnical Journal}, 58\penalty0 (10):\penalty0 1466--1477, 2021.

\bibitem[Robinson et~al.(2019)Robinson, Brown, Matsui, Brennan, Augarde, Coombs, and Cortis]{robinson2019centrifuge}
Scott Robinson, Michael~John Brown, Hidetake Matsui, Andrew Brennan, Charles Augarde, Will Coombs, and Michael Cortis.
\newblock Centrifuge testing to verify scaling of offshore pipeline ploughs.
\newblock \emph{International Journal of Physical Modelling in Geotechnics}, 19\penalty0 (6):\penalty0 305--317, 2019.

\bibitem[Schanz et~al.(2019)Schanz, Vermeer, and Bonnier]{schanz2019hardening}
T.~Schanz, P.A. Vermeer, and P.G. Bonnier.
\newblock The hardening soil model: Formulation and verification.
\newblock In \emph{Beyond 2000 in computational geotechnics}, pages 281--296. Routledge, 2019.

\bibitem[Simo(1992)]{Simo1992a}
JC~Simo.
\newblock Algorithms for static and dynamic multiplicative plasticity that preserve the classical return mapping schemes of the infinitesimal theory.
\newblock \emph{Computer Methods in Applied Mechanics and Engineering}, 99:\penalty0 61--112, 1992.

\bibitem[Singer et~al.(2024)Singer, Teschemacher, Larese, W\"{u}chner, and Bletzinger]{Singer2024}
V.~Singer, T.~Teschemacher, A.~Larese, R.~W\"{u}chner, and K.-U. Bletzinger.
\newblock Lagrange multiplier imposition of non-conforming essential boundary conditions in implicit material point method.
\newblock \emph{Computational Mechanics}, 73:\penalty0 1311–1333, 2024.

\bibitem[Solowski et~al.(2021)Solowski, Berzins, Coombs, Guilkey, Möller, Tran, Adibaskoro, Seyedan, Tielen, and Soga]{Solowski2021}
Wojciech~T. Solowski, Martin Berzins, William~M. Coombs, James~E. Guilkey, Matthias Möller, Quoc~Anh Tran, Tito Adibaskoro, Seyedmohammadjavad Seyedan, Roel Tielen, and Kenichi Soga.
\newblock Material point method: Overview and challenges ahead.
\newblock volume~54 of \emph{Advances in Applied Mechanics}, pages 113--204. Elsevier, 2021.

\bibitem[Sticko et~al.(2020)Sticko, Ludvigsson, and Kreiss]{Sticko2020}
S.~Sticko, G.~Ludvigsson, and G.~Kreiss.
\newblock High-order cut finite elements for the elastic wave equation.
\newblock \emph{Advances in Computational Mathematics}, 46:45:\penalty0 1 -- 28, 2020.

\bibitem[Sulsky et~al.(1994)Sulsky, Chen, and Schreyer]{sulsky1994particle}
D.~Sulsky, Z.~Chen, and H.L. Schreyer.
\newblock A particle method for history-dependent materials.
\newblock \emph{Computer Methods in Applied Mechanics and Engineering}, 118\penalty0 (1-2):\penalty0 179--196, 1994.

\bibitem[Sulsky et~al.(1995)Sulsky, Zhou, and Schreyer]{Sulsky1995}
Deborah Sulsky, Shi-Jian Zhou, and Howard~L Schreyer.
\newblock Application of a particle-in-cell method to solid mechanics.
\newblock \emph{Computer Physics Communications}, 87\penalty0 (1):\penalty0 236--252, 1995.

\bibitem[Tielen et~al.(2017)Tielen, Wobbes, Möller, and Beuth]{Tielen2017}
Roel Tielen, Elizaveta Wobbes, Matthias Möller, and Lars Beuth.
\newblock A high order material point method.
\newblock \emph{Procedia Engineering}, 175:\penalty0 265--272, 2017.
\newblock Proceedings of the 1st International Conference on the Material Point Method (MPM 2017).

\bibitem[Wang et~al.(2019)Wang, Coombs, Augarde, Cortis, Charlton, Brown, Knappett, Brennan, Davidson, Richards, and Blake]{Wang2019}
L.~Wang, W.M. Coombs, C.E. Augarde, M.~Cortis, T.J. Charlton, M.J. Brown, J.~Knappett, A.~Brennan, C.~Davidson, D.~Richards, and A.~Blake.
\newblock On the use of domain-based material point methods for problems involving large distortion.
\newblock \emph{Computer Methods in Applied Mechanics and Engineering}, 355:\penalty0 1003--1025, 2019.

\bibitem[Wang et~al.(2021)Wang, Coombs, Augarde, Cortis, Brown, Brennan, Knappett, Davidson, Richards, White, et~al.]{wang2021efficient}
Lei Wang, William~M Coombs, Charles~E Augarde, Michael Cortis, Michael~J Brown, Andrew~J Brennan, Jonathan~A Knappett, Craig Davidson, David Richards, David~J White, et~al.
\newblock An efficient and locking-free material point method for three-dimensional analysis with simplex elements.
\newblock \emph{International Journal for Numerical Methods in Engineering}, 122\penalty0 (15):\penalty0 3876--3899, 2021.

\bibitem[Wriggers(2006)]{wriggers2006computational}
P.~Wriggers.
\newblock \emph{Computational contact mechanics}.
\newblock Springer, 2006.

\bibitem[Yamaguchi et~al.(2021)Yamaguchi, Moriguchi, and Terada]{Yamaguchi2021}
Yuya Yamaguchi, Shuji Moriguchi, and Kenjiro Terada.
\newblock Extended b-spline-based implicit material point method.
\newblock \emph{International Journal for Numerical Methods in Engineering}, 122\penalty0 (7):\penalty0 1746--1769, 2021.

\end{thebibliography}

\appendix
\section{Rigid body first and second variations}\label{App:truss frame}
This appendix provides the first and second variations of a point on a triangle surface, $\bm{x}^\prime$. The definition from \eqref{equ:app xn definition m} is repeated here, with a slight abuse of the sum notation, 
\begin{equation}\label{equ:app Nx full sum}
    \bm{x}^\prime = N(\xi)\bm{G}(\bm{\theta}) = (A \bm{R} + B\bm{I}):\bm{t}_{RB} + \bm{x}_{M}.
\end{equation}
Its first variation is 
\begin{equation}\label{equ:dx first variation}
    \delta\bm{x}^\prime = \frac{\partial N(\xi)}{\partial \xi^\alpha}\delta\xi^\alpha\bm{G}(\bm{\theta})
    +  N(\xi)\frac{\partial\bm{G}(\bm{\theta})}{\partial \bm{\theta}}\cdot\delta\bm{\theta},
\end{equation}
and second variation
\begin{equation}\label{equ:dx second variation}
    \Delta\delta\bm{x}^\prime = 
    \underbrace{\frac{\partial N(\xi)}{\partial \xi^\alpha}\Delta\delta\xi^\alpha\bm{G}(\bm{\theta})}_{:=\bm{t}_\alpha\Delta\delta\xi^\alpha}
    + \underbrace{\frac{\partial N(\xi)}{\partial \xi^\alpha}\delta\xi^\alpha\frac{\partial\bm{G}(\bm{\theta})}{\partial\bm{\theta}}\Delta\theta}_{:= \left(\frac{\partial(\Delta\bm{x}^\prime)}{\partial \xi^\beta}\right)\delta\xi^\alpha}
    + \underbrace{\frac{\partial N(\xi)}{\partial \xi^\alpha}\Delta\xi^\alpha\frac{\partial\bm{G}(\bm{\theta})}{\partial\bm{\theta}}\delta\theta}_{:= \left(\frac{\partial(\delta\bm{x}^\prime)}{\partial \xi^\beta}\right)\Delta\xi^\alpha}
    +   \underbrace{N(\xi)\Delta\bm{\theta}\cdot\frac{\partial^2\bm{G}(\bm{\theta})}{\partial \bm{\theta}^2}\cdot\delta\bm{\theta}}_{:= \delta\bm{\theta}\cdot\left(\frac{\partial^2\bm{x}^\prime}{\partial\bm{\theta}^2}\right)\cdot\Delta\bm{\theta}}.
\end{equation}
which requires the following definitions,
\begin{equation}\label{equ:div dx first variation}
    \frac{\partial(\delta\bm{x}^\prime)}{\partial \xi^\beta} =  \delta\left(\frac{\partial\bm{x}^\prime}{\partial \xi^\beta}\right)
    \quad\text{and}\quad
   \frac{\partial(\Delta\bm{x}^\prime)}{\partial \xi^\beta} =  \Delta\left(\frac{\partial\bm{x}^\prime}{\partial \xi^\beta}\right)
\end{equation}
and
\begin{equation}
    \Delta\left(\frac{\partial(\delta\bm{x}^\prime)}{\partial \xi^\beta}\right) = \frac{\partial N(\xi)}{\partial \xi_\beta}\Delta\bm{\theta}\cdot\frac{\partial^2\bm{G}(\bm{\theta})}{\partial \bm{\theta}^2}\cdot\delta\bm{\theta}.
\end{equation}
Additionally Equations \eqref{equ:dx first variation} and \eqref{equ:dx second variation} are dependent on the first and second variations of $\bm{G}(\bm{\theta})$, respectively given as
\begin{equation}
    \frac{\partial\bm{G}_n(\bm{\theta})}{\partial\bm{\theta}}\cdot\delta\bm{\theta} = (A_n\bm{R} + B_n\bm{I}):\delta\bm{t}_{RB} + \delta\bm{x}_{M}
\end{equation}
and
\begin{equation}
    \Delta\bm{\theta}\cdot\frac{\partial^2\bm{G}_n(\bm{\theta})}{\partial\bm{\theta}^2}\cdot\delta\bm{\theta} = (A_n\bm{R} + B_n\bm{I}):\Delta\delta\bm{t}_{RB}.
\end{equation}
Furthermore the first and second variation of $\bm{t}_{RB}$ are also required, respectively,
\begin{equation}
\delta\bm{t}_{RB} = \left(\frac{\bm{I} - \bm{t}_{RB}\otimes\bm{t}_{RB}}{\|\bm{v}\|}\right)\cdot\frac{\partial \bm{v}}{\partial t}\delta t
    = \left(\frac{\bm{I} - \bm{t}_{RB}\otimes\bm{t}_{RB}}{\|\bm{v}\|}\right)\cdot\delta \bm{v}
    \quad\text{where}\quad\bm{v} = \bm{x}_{D} - \bm{x}_{M}.
\end{equation}
and
\begin{equation}
    \begin{split}
         \Delta\delta\bm{t}_{RB} &= \Delta\left[\left(\frac{\bm{I} - \bm{t}_{RB}\otimes\bm{t}_{RB}}{\|\bm{v}\|}\right)\cdot\delta \bm{v}\right]\\
         & = \Delta\left(\frac{1}{\|\bm{v}\|}\right)\left(\bm{I} - \bm{t}_{RB}\otimes\bm{t}_{RB}\right)\cdot\delta \bm{v}
           + \left(\frac{1}{\|\bm{v}\|}\right)\left(\bm{I} - \Delta\bm{t}_{RB}\otimes\bm{t}_{RB}\right)\cdot\delta \bm{v}
           + \left(\frac{1}{\|\bm{v}\|}\right)\left(\bm{I} - \bm{t}_{RB}\otimes\Delta\bm{t}_{RB}\right)\cdot\delta \bm{v},
    \end{split}
\end{equation}
where
\begin{equation*}
    \Delta\left(\frac{1}{\|\bm{v}\|}\right) = \frac{-\bm{v}}{\|\bm{v}\|^3}\cdot\Delta\bm{v}\quad\text{and}\quad\Delta \bm{t}_{RB} = \left(\frac{\bm{I} - \bm{t}_{RB}\otimes\bm{t}_{RB}}{\|\bm{v}\|}\right)\cdot\Delta\bm{v}.
\end{equation*}

\section{Contact first and second variations}

This appendix provides the first and second variations of the contact components required to implement the proposed contact method.  

\subsection{Gap function}\label{App:gap var}
Following closely the work of Pietrzak and Curnier \cite{pietrzak1999large}, and Wriggers \cite{wriggers2006computational}, the forms of the first and second variations of the gap functions can be provided. The gap function, \eqref{equ:gap function}, is restated here for convenience
\begin{equation}
    g_N = \left[\bm{x}(\tau) - \bm{x}^\prime(\xi(\tau)^\alpha,\tau)\right]\cdot\bm{n}
\end{equation}
where $\bm{x}$ is the point, the vertex of the GIMP domain, in contact with the rigid body surface at time $\tau$, where $\xi^\alpha(\tau)$, $\alpha\in\{1,2\}$, is a function describing the local position of the CPP projection onto the surface that is used to describe the global contact position $\bm{x}^\prime$.

The first variation of the gap function is
\begin{equation}
    \delta g_N = \left[\delta\bm{x} - \delta\bm{x}^\prime \right]\cdot \bm{n}
\end{equation}
where $\delta$ denotes the first variation. The second variation of the normal gap $g_N$ is
\begin{equation}
    \Delta\delta g_n = -\bm{n}\cdot\left(\frac{\partial(\delta \bm{x}^\prime)}{\partial \xi^\alpha}\Delta\xi^\alpha + \frac{\partial(\Delta \bm{x}^\prime)}{\partial \xi^\alpha}\delta\xi^\alpha+ \delta\bm{\theta}\cdot\left(\frac{\partial^2\bm{x}^\prime}{\partial\bm{\theta}^2}\right)\cdot\Delta\bm{\theta}\right)
    +
    g_N\bm{n}\cdot\left(\frac{\partial(\delta \bm{x}^\prime)}{\partial \xi^\alpha}\right)A^{\alpha\beta}\left(\frac{\partial(\Delta \bm{x}^\prime)}{\partial \xi^\beta}\right)\cdot\bm{n}
\end{equation}
where $\xi^\alpha$ is the local coordinate on triangular element of the rigid body, with its first variations $\delta\xi$ or $\Delta\xi$, and 
\[[A^{\alpha\beta}]^{-1} = [A_{\alpha\beta}] = \bm{t}_\alpha\cdot\bm{t}_\beta\]
is the inverse of the first fundamental form matrix.

\subsection{Variations of $\xi$}\label{App:xi var}
The first and second variations of the function describing the local coordinate $\xi^\alpha(\tau)$ of the CPP are also required. The first variation is described as 
\begin{equation}
    \delta\xi^\alpha = A^{\alpha\beta}\left[(\delta\bm{x} - \delta\bm{x}^\prime)\cdot\bm{t}_\beta + g_N\bm{n}\cdot\frac{\partial(\delta \bm{x}^\prime)}{\partial\xi^\beta}  \right]
\end{equation}
and the corresponding second variation is
\begin{equation}
\begin{split}
    \Delta\delta\xi^\alpha = A^{\alpha\beta}\biggl[&
    -\bm{t}\cdot\left(
    \delta\xi^\sigma\frac{\partial(\Delta \bm{x}^\prime)}{\partial\xi^\sigma}
    + \Delta\xi^\gamma\frac{\partial(\delta \bm{x}^\prime)}{\partial\xi^\gamma}    
    + \delta\bm{\theta}\left(\frac{\partial^2\bm{x}^\prime}{\partial\bm{\theta}^2}\right)\Delta\bm{\theta}
    \right)\\
    &-\delta\xi^\gamma\bm{t}_{\gamma}
    \left(
    \frac{\partial(\Delta \bm{x}^\prime)}{\partial\xi^\beta}
    \right)
     -\Delta\xi^\gamma\bm{t}_{\gamma}
    \left(
     \frac{\partial(\delta \bm{x}^\prime)}{\partial\xi^\beta}
    \right)\\
    &+(\delta\bm{x} - \delta\bm{x}^\prime) \frac{\partial(\Delta \bm{x}^\prime)}{\partial\xi^\beta} 
    +(\Delta\bm{x} - \Delta\bm{x}^\prime) \frac{\partial(\delta \bm{x}^\prime)}{\partial\xi^\beta}\\
    & + g_N\bm{n}\cdot\Delta\left(\frac{\partial(\delta {\bm{x}^\prime})}{\partial\xi^\beta}\right)\biggr].
\end{split}
\end{equation}

\subsection{Stick and slip}\label{App:stick and slip}
The first variation of the frictional load in the tangential direction, $\bm{t}_{\alpha}\cdot\bm{p}_T$, needs to be determined for the stick and slip states. Starting generally, the first variation is 
\begin{equation}\label{equ:app dpt}
    \Delta\left(\bm{t}_{\alpha}\cdot\bm{p}_T\right) = \left[\left(\Delta\bm{t}_{\alpha}\right)\cdot\bm{p}_T\right]
    +
    \left[\bm{t}_{\alpha}\cdot\left(\Delta\bm{p}_T\right)\right].
\end{equation}
The $\Delta t_\alpha$ term is simply,
\begin{equation}
    \Delta t_\alpha = \frac{\partial(\Delta \bm{x}^\prime)}{\partial\xi^\alpha}
\end{equation}
noting that $\delta$ and $\Delta$ are equivalent operators. The second term of \eqref{equ:app dpt} requires the variation $\Delta\bm{p}_T$, which is either acting in stick or slip, defined in Equation \eqref{equ:stick slip}.

\subsubsection{Stick}
For the \textbf{stick} case the frictional force takes the form 
\begin{equation*}
    \bm{p}_T = \bm{p}_{T,m} + \epsilon_T\Delta\bm{g}_T =  \bm{p}_{T,m} + \epsilon_T\bm{g}_T^\Delta
\end{equation*}
where $\Delta\bm{g}_T$, the total tangential movement over the surface during a time step, is redefined as $\bm{g}_T^\Delta$ to avoid conflicting notation with first and second variations. It is first convenient to define the variation in tangential movement,
\begin{equation}
    \Delta\bm{g}_T^\Delta = (\Delta\bm{t}_{\alpha,m+1})(\xi_{m+1}^\alpha + \xi_{m}^\alpha) + \bm{t}_{\alpha,m+1}(\Delta\xi_{m+1}^\alpha).
\end{equation}
This enables the definition for the first variation of the stick contact force,
\begin{equation}
    \Delta\bm{p}_T = \epsilon_T\Delta\bm{g}_T^\Delta
\end{equation}

\subsubsection{Slip}
For the \textbf{slip} case $\bm{p}_T$ has the form
\begin{equation*}
    \bm{p}_T = \mu|p_N|\left(\frac{{\bm{p}_{tr}}}{||{\bm{p}_{tr}}||}\right)
\end{equation*}
and its corresponding variation is
\begin{equation}
\begin{split}
    \Delta\bm{p}_T &= \mu(\Delta|p_N|)\left(\frac{{\bm{p}_{tr}}}{||{\bm{p}_{tr}}||}\right) +  \mu|p_N|\Delta\left(\frac{{\bm{p}_{tr}}}{||{\bm{p}_{tr}}||}\right)\\
\end{split}
\end{equation}
where
\begin{equation}
    \Delta|p_N| = \epsilon_N\delta|g_N| = \epsilon_N\frac{g_N}{|g_N|}\cdot\Delta g_N
\end{equation}
and
\begin{equation}
    \Delta\left(\frac{{\bm{p}_{tr}}}{||{\bm{p}_{tr}}||}\right) = \left(\frac{\bm{I}-\bm{s}\otimes\bm{s}}{\|\bm{p}_{tr}\|}\right)\cdot\Delta \bm{p}_{tr}\quad\text{where}\quad\bm{s}=\bm{p}_{tr}/\|\bm{p}_{tr}\|
\end{equation}
noting that $\bm{p}_{tr} =  \bm{p}_{T,m} + \epsilon_T\bm{g}_T^\Delta$; it is the same as the stick form of $\bm{p}_{T}$ and also therefore has the same form for the first variation.

\end{document}